%% file: main.tex
\documentclass{article} 
\usepackage{iclr2027_conference,times}
\usepackage{graphicx}
\usepackage{booktabs}
\graphicspath{ {figures/} }

\input{math_commands.tex}

\usepackage{hyperref}
\usepackage{url}
\usepackage{listings}
\usepackage{titletoc}
\usepackage{xcolor}

\newif\ifextended

\title{Discovering Practical Algorithms for \\ Least Squares using Language Models}

\author{Marco Pirazzini\thanks{Corresponding author. Work performed during an internship at Sakana AI.} \\
Yale University\\
\texttt{marco.pirazzini@yale.edu} \\
\And
Robert Tjarko Lange \\
Sakana AI \\
\texttt{robert@sakana.ai} 
}

\usepackage{cleveref}

\iclrfinalcopy 
\renewcommand{\headrulewidth}{0pt}

\begin{document}

\maketitle

\fancyhead{}

\begin{abstract}
Least squares is a central problem in scientific computing with applications across many domains such as the natural and social sciences. The design and implementation of state-of-the-art algorithms for solving least squares are often based on worst-case theoretical analysis and tailored to provide numerically stable solutions. These principles are often too conservative in real-world settings, incurring severe runtime slowdowns. Manually designing per-use-case algorithm implementations, on the other hand, requires domain knowledge and tedious trial-and-error. In this work, we leverage recent progress in large language model-guided automated discovery to automate the design of practical algorithms for numerical linear algebra applications. More specifically, our methodology frames the search for practical least squares algorithms as a meta-discovery problem: given a family of downstream solver tasks, we combine language model-guided algorithm evolution with robust meta-evaluation accounting for runtime and solver accuracy. We demonstrate the effectiveness and flexibility of our framework by discovering faster practical algorithms across an extensive list of synthetic and real-world least squares tasks. The discovered algorithms are $2\times$ faster than optimized baselines on large kernel regression tasks from the SUSY dataset for high-energy particle physics classification with up to one million rows. Additionally, we analyze the algorithms discovered by our framework and decompose the runtime improvements across algorithmic changes, hyperparameter optimization, and code engineering.
\end{abstract}

\section{Introduction}\label{sec::intro}

Data-driven fitting of a linear relationship by minimizing the squared errors is one of the most studied problems in scientific computing, with wide-ranging applications including machine learning \citep{hastie2009elements, rudi2017falkon}, economics \citep{wooldridge2009introductory}, physics \citep{bevington2003data, aster2005parameter}, medical sciences \citep{hansen2021computed, vittinghoff2004regression}, astronomy \citep{feigelson2012modern} and many more. This is known as the \emph{least squares problem} (LS), where given a matrix $A \in \mathbb{R}^{m \times n}$ with $m > n$ and a vector $b \in \mathbb{R}^m$, the goal is to find
\begin{equation}\label{eq::leastsquares}
    x^*=\arg\min_{x \in \mathbb{R}^n} \|Ax - b\|_2^2.
\end{equation}
Many algorithms with strong theoretical guarantees have been developed to solve this problem, such as backward-stable direct solvers based on QR or SVD factorization \citep{golub1965numerical}, faster but less accurate direct solvers based on the Cholesky factorization of the Gram matrix \citep{golub2013matrix}, and randomized preconditioned iterative algorithms based on sketching \citep{murray2023randnla}. They have highly optimized software implementations in C/Fortran that are directly wrapped and made available in popular programming languages \citep{dongarra1990set, anderson1999lapack}.

\begin{figure}[h!]
\begin{center}
    \includegraphics[width=0.9\textwidth]{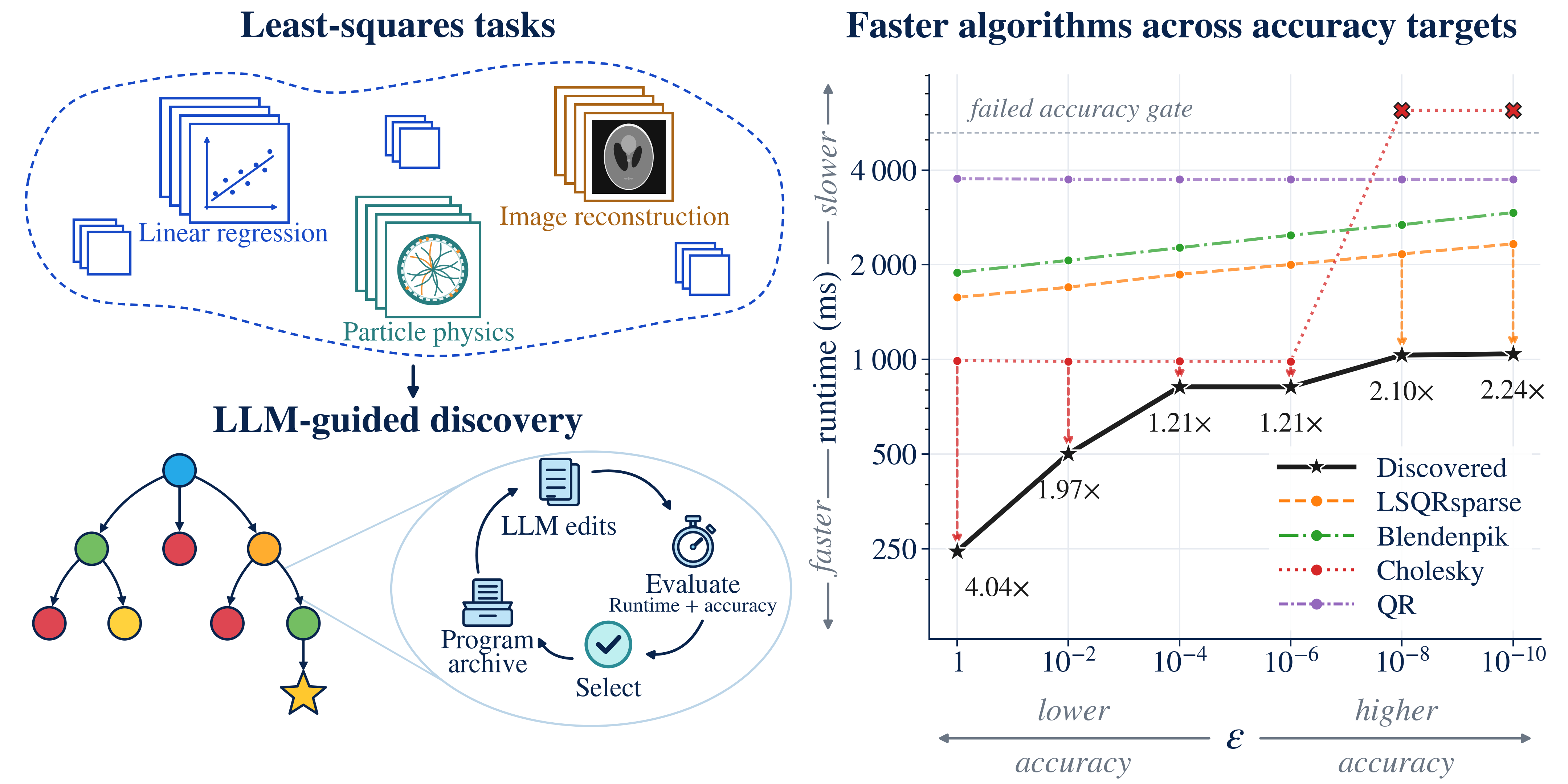}
\end{center}
\caption{
\textbf{Left.} Given a specific family of least squares problem classes, we run LLM-guided evolution to discover faster practical algorithms tailored to the task. \textbf{Right.} The method discovers algorithms faster than established baselines on a variety of tasks, for example, where we vary the target accuracy level for large ill-conditioned matrices with size $(20000,2000)$ and conditioning $10^6$. The best baseline changes with the task, but our framework adapts to the task and can discover tailored improvements that are especially large in low-accuracy regimes.}
\label{fig::1}
\end{figure}

Despite the maturity of this field from both a mathematical and software perspective, the gap between theory and practice is wide when it comes to runtime speed (\Cref{fig::1}, right). Many algorithms and implementations are developed with a worst-case perspective, where the priority is to find numerically stable high-accuracy solutions on all inputs. This leads to conservative choices that are often excessively slow for practical applications. The worst-case viewpoint is ill-suited to bridge this gap: the space of feasible LS solvers is large and includes a variety of algorithmic families with different trade-offs. Overly conservative worst-case settings often neglect better algorithmic designs, which outperform on (still broad) application domains. 

Human researchers are limited by time and domain knowledge and hence often do not spend suitable amounts of time optimizing the tools used. Given the advent of automated LLM-guided algorithm discovery, we argue that investing some time in searching for a better task-specific algorithm is now more worthwhile, and can lead to significant runtime improvement for downstream solves. Here, we frame the problem of designing dedicated LS algorithms as a meta-discovery problem (\Cref{fig::1}): given a process for generating task instances, we leverage offline test-time scaling with verifier feedback information to discover suitable high-performance algorithms. 

We demonstrate that this approach is capable of finding runtime improvements across a variety of accuracy levels and across different problem instance families. The cost of the discovery process can later be amortized by using the resulting optimized algorithm to solve the downstream task repeatedly.
We thereby move from a worst-case to a \textit{practical-case} perspective on algorithm design for least squares. In a \textit{practical least squares task}, a user has a specific workload that generates instances from a fixed data-generating process. Average-case analysis in which one analyzes an algorithm on random inputs avoids some pitfalls of worst-case analysis, but can lead to poor algorithms for inputs generated by different processes. Therefore, the design of truly practical algorithms must be data-driven, where the resulting algorithm is adapted to instances from the data-generating process of the given task.
In summary, we leverage the recent developments in LLM-guided evolutionary search to discover faster practical algorithms tailored to specific applications (\Cref{sec::shinkalgodisc}). The main contributions of our work are the following:
\begin{itemize}
    \item \textbf{Worst-case versus practical-case LS algorithm design framework} (\Cref{sec::tasks}): We introduce an applied perspective to algorithm design and demonstrate how to illuminate the space of least squares algorithms. To the best of our knowledge, this provides the first systematic study of LLM-guided algorithm evolution applied to practical LS problems.

    \item \textbf{Runtime improvement discoveries} (\Cref{sec::experiments}): Our method discovers faster algorithms for a wide array of synthetic LS problems, real-world tasks and target solution accuracy levels. The efficiency gains are especially strong for large-scale or low-accuracy settings.

    \item \textbf{Generalization \& discovery analysis} (\Cref{sec::analysis}): Our method can be flexibly applied to diverse settings and requires minimal work to customize. Additionally, we demonstrate that the discovered algorithms generalize beyond the meta-training task distribution. Finally, we deconstruct algorithmic improvements and contrast LLM-guided discovery with black-box hyperparameter tuning.
\end{itemize}

\section{Background and Related Work}\label{sec::background}

\paragraph{Least Squares Solvers.}
Least squares is a well-studied computational problem with several algorithmic families developed to find good solutions. A standard robust solver for (\ref{eq::leastsquares}) when the matrix has full column rank is based on a Householder QR decomposition followed by a triangular solve \citep{golub1965numerical}:
\begin{align*}
    A& = QR \\
    x^* &= R^{-1}Q^\top b
\end{align*}
Other direct methods are based on pivoted QR, SVD, or Cholesky on the Gram matrix $A^\top A$ (see \citet{bjorck2024numerical} for a modern treatment). Sketching can be used to speed up direct methods or preconditioned iterative methods such as conjugate gradient on the normal equations \citep[CGLS,][]{hestenes1952methods} or LSQR \citep{paige1982lsqr}. The core idea consists of drawing a random matrix $S \in \mathbb{R}^{d \times m}$ with $n \le d < m$, and compressing the information of $A$ into the smaller matrix $A_s = SA$ \citep{murray2023randnla,martinsson2020randomized}. A simplified description of a sketch-and-precondition algorithm running for $k$ iterations with a QR-based preconditioner is the following.
\begin{align*}
    A_s &= Q_s R_s \\
    y_0 &= Q_s^\top S b \implies x_0 = R_s^{-1} y_0 \quad \text{(sketch-and-solve solution)}\\
    y_k &= \textsc{LSQR}\left( AR_s^{-1}, b; y_0, k\right) \\
    x_k &= R_s^{-1}y_k \quad \text{(sketch-and-precondition solution)}
\end{align*}
Early work on randomized preconditioners by \citet{rokhlin2008fast} led to faster practical algorithms \citep{avron2010blendenpik, meng2014lsrn, cartis2021hashing}. More recent works on adaptive preconditioners \citep{huh2026adaptive} and backward-stable randomized solvers \citep{epperly2026fast} recommend sparse sketching and warm-starting as defaults for fast algorithms.

\paragraph{LLM-Guided Evolutionary Program Discovery.}\label{sec::llm-background}
LLM-guided evolutionary search as a tool for algorithm discovery offers an alternative avenue to traditional genetic programming or Reinforcement Learning for optimizing over the space of algorithms by representing them as computer programs and leveraging LLMs as mutation operators. This is a very active research area popularized by FunSearch \citep{romera2024mathematical} and refined by many works including AlphaEvolve \citep{novikov2025alphaevolve}, OpenEvolve \citep{sharma2025openevolve}, ShinkaEvolve \citep{lange2026shinkaevolve}, and SkyDiscover \citep{liu2026skydiscover}. Even more recently, \citet{aoudia2026ai} extended these ideas by using coding agents instead of individual API calls. Most applications of LLM-guided search to mathematics and theoretical computer science focus on evaluating the output of the algorithms (e.g., as a certificate for some inequality), rather than properties of the algorithm itself; see \cite{georgiev2025mathematical}. AlphaEvolve \citep{novikov2025alphaevolve} and AlgoTuner \citep{press2026algotune} study tasks where the goal is to optimize the algorithm's runtime, but neither considers general least squares problems.
LLM-guided evolutionary algorithms have also been studied in the context of Kalman filtering \citep{greenberg2023optimization, saketos2026data}.

\paragraph{Algorithm Discovery for Least Squares.}
A closely related work on data-driven discovery of practical algorithms for least squares is \cite{cho2025surrogate}. Their approach leverages Gaussian Process-based black-box tuning to find favorable parameter configurations for a certain family of sketch-and-precondition algorithms. Their approach requires specifying the hyperparameters, categorical variables, and admissible values a priori, and is only studied in the context of randomized iterative algorithms. Generalizing this to fully capture the space of algorithms for least squares and the available hardware/software details is laborious and very difficult to achieve. Our approach provides a flexible alternative: algorithms are represented as computer programs and evolved using frontier LLMs. 

More recently, \cite{xiong2026rl4rla} used reinforcement learning to discover symbolic programs from a specified curriculum of increasingly difficult tasks. Their method succeeds in re-discovering known high-level algorithms, which highlights the interpretability of symbolic program discovery. Our settings are different in that we do not specify the prescribed sequence of algorithmic properties via a curriculum, and we use LLM-guided evolution to find a faster practical algorithm tailored to the task of interest. We will see that the algorithms discovered by our framework are also interpretable.

\section{Task \& Discovery Methodology}\label{sec::shinkalgodisc}

The main methodological contribution of our work consists in designing a simple framework for evaluating the efficiency of algorithms on a specific task of interest, provided such task has a way to verify that the output of the algorithm satisfies the target. This applies well beyond least squares, but we focus on this application for its practical relevance. In Section \ref{sec::experiments}, we provide evidence that our method works well with minimal changes across a wide range of least squares problems.

\paragraph{Least Squares Algorithm Discovery Tasks.}\label{sec::tasks}
We define a \emph{least squares task} as the tuple $\mathcal{T} = (I_{\text{train}}, I_{\text{test}}, \varepsilon, \texttt{msg})$, where $I_{\text{train}}=\{((A_1,b_1),x^*_1), ..., ((A_K,b_K), x^*_K)\}$ and $I_{\text{test}}$ are sets of least squares data instances, $\varepsilon > 0$ is the accuracy target and \texttt{msg} is a natural language description of the task. $I_{\text{train}}$ is used during the meta-discovery phase to ``train'' an algorithm for $\mathcal{T}$, while $I_{\text{test}}$ acts as a set of held-out instance-solution pairs, which are used to study the generalization quality of the discovered algorithm.
Throughout this paper, we focus on dense, full-rank and inconsistent (i.e. $\|Ax^*-b\| > 0$) least squares problems.  There are different notions of accuracy for a solution to a least squares problem as in (\ref{eq::leastsquares}). In this paper, we follow the standard $(1+\varepsilon)$-approximation criterion in randomized numerical linear algebra \citep{clarkson2017low} and consider the relative residual ratio:
\begin{equation}\label{eq::acc_metric}
    f_{\text{res}}(x) = \frac{\|Ax - b\|_2}{\|Ax^* - b\|_2}.
\end{equation}
We say that an algorithm $\mathcal{A}$ is feasible for solving the task $\mathcal{T}$ if the output of $\mathcal{A}$ on all instances of $\mathcal{T}$ satisfies $f_{\text{res}} \le 1+\varepsilon$. 
Other important notions of accuracy for least squares algorithms are forward error and backward error, and many of the baselines discussed above are designed to minimize these errors. Modifying the notion of feasibility according to these metrics is straightforward.

\paragraph{LLM-Guided Program Evolution Framework.}
We explore the space of algorithms for a given least squares task $\mathcal{T}$, which is a large and non-differentiable search space. Two components are needed to successfully optimize over this space: a way to effectively explore it by generating candidate algorithms \textbf{(C1)}, and a way to accurately and robustly score a candidate algorithm \textbf{(C2)}.

Traditionally, this search was done manually by human hypothesis testing. Here, we aim to automate the design of suitable algorithms based on automated iterative program refinement. Due to their strong coding and mathematical reasoning capabilities, modern LLMs offer a very interesting black-box sampling oracle for component \textbf{(C1)}. The first natural thing to try is to just sample many times from (an ensemble of) LLMs and keep the best program (\textit{best-of-N}).
A different approach is to use test-time scaling frameworks inspired by evolutionary computation to balance exploration with exploitation of different ideas already attempted. We adopt the latter \citep{lange2026shinkaevolve} but validate other alternatives (see Appendix \ref{app::shinka_vs_restarts}).
The initial program used for every task presented in this paper was always the same: a randomized preconditioned iterative algorithm using a sparse sketch to build the preconditioner followed by LSQR. The task description ``\texttt{msg}'' was essentially fixed, with minimal changes describing the desired accuracy target, see Appendix \ref{app::method} for more details.

The natural candidate for component \textbf{(C2)} is the typical runtime of the algorithm on the target matrices, conditional on the output of the algorithm achieving sufficient accuracy. This poses several challenges, because measuring program runtime directly often gives a very noisy estimate of true runtime, due to the state of the machine at the moment of running (e.g. contention, thermal state, cold vs warm cache, etc.).

\paragraph{The Least-Squares Discovery Objective Function.}
Let $\mathcal{A}$ be the algorithm currently being evaluated, $\mathcal{A}_0$ be the initial algorithm, and $\mathcal{T}$ be the task of interest as defined in Section \ref{sec::tasks}. For training matrices $(A_i,b_i)$ in $\mathcal{T}$, we run algorithms $\mathcal{A}$ and $\mathcal{A}_0$ $R$ times, and let $t_\mathcal{A}(i,r)$ be the runtime of $\mathcal{A}(A_i,b_i)$ on run $r$. Let $\rho(i,r) = \frac{t_\mathcal{A}(i,r)}{t_{\mathcal{A}_0}(i,r)}$ be the ratio of the runtime of the current program over the starting one, then we define the \emph{log geometric-mean speedup} of $\mathcal{A}$ as
\begin{equation}\label{eq::shinka-score-corr}
    f(\mathcal{A}) = -\frac{1}{K}\sum_{i=1}^K \log_2 \left(\text{median}_{r=1,...,R}(\rho(i,r)) \right)
\end{equation}

Focusing on the ratio $\rho$ relative to a fixed algorithm instead of raw runtime $t$ has two benefits. First, the normalization provides a scale-invariant objective, which gives each instance equal weight in terms of multiplicative speedup over $\mathcal{A}_0$. Second, re-timing the initial program each time lowers the fluctuations due to machine contention at each generation. Similarly, running an algorithm $R$ times has two benefits. First, many practical least squares algorithms are randomized, so we can increase our confidence in the algorithm's correctness by running multiple trials. Second, taking the median among restarts reduces the timing noise.\footnote{The minimum is a popular choice of aggregator for dealing with runtime noise, because contention can only make programs slower. We choose the median because some algorithms are randomized and it removes outliers and avoids selecting lucky seeds.} More details on other measures taken to reduce runtime noise are in Appendix \ref{app::method}.

Some candidate algorithms might fail to reach the target accuracy $\varepsilon$ on all matrices and all runs. A simple way to deal with incorrect programs is to assign score $-\infty$, but this is too non-smooth and does not allow close-to-working ideas to survive. We say that an instance $(A_i, b_i)$ is ``solved'' by $\mathcal{A}$ if all $R$ runs reach the desired accuracy $f_{\text{res}} \le 1+\varepsilon$, and ``unsolved'' otherwise, and we add a soft penalty that gracefully degrades with the number of unsolved instances. Let $S$ be the subset of solved instances, $\tau < 0$ be a parameter demarcating the threshold between correct and incorrect algorithms, and $g(\mathcal{A},S)>0$ be a penalty function defined in Appendix \ref{app::method}. Consider the score function
\begin{equation}\label{eq::shinka-score}
    \text{score}(\mathcal{A}) = \begin{cases}
        \max\{f(\mathcal{A}), \tau\} & \text{if } |S| = K\\
        \tau - g(\mathcal{A},S) & \text{if } |S| < K
    \end{cases}
\end{equation}
Our discovery framework combines \textbf{(C1)} and \textbf{(C2)} by using ShinkaEvolve \citep{lange2026shinkaevolve} to generate algorithms $\mathcal{A}$ with evaluation function (\ref{eq::shinka-score}). The only information from $\mathcal{T}$ that is exposed to the LLM are \texttt{msg}, $\varepsilon$ and some aggregate information from $I_{\text{train}}$ such as average instance size or conditioning, while the full instances and solutions inside $I_{\text{train}}$ are used privately by the evaluator to compute the score (\ref{eq::shinka-score}). $I_{\text{test}}$ is only used for held-out testing.

The total computational cost of the discovery framework consists of $2RK$ least squares solves per candidate evaluation, excluding the cost of collecting the instances in $I_{\text{train}}$ and computing their reference solutions. The exact time required varies depending on the configuration chosen and the complexity of the task considered, and the amount of time waiting for a response from the LLM queries. For example, large-scale, high-accuracy tasks take longer to be solved so they make the whole discovery pipeline slower. Moreover, if someone has a machine with very controlled minimal contention, they can speed up evaluation by removing co-timing or reducing the number of repeats. As a representative, the tasks discussed in Figure \ref{fig::4.2} below required on average around $3$ hours to produce the final algorithm, using around $\$10$ in API spend. The dollar cost was consistently around that amount on all experiments discussed below, because it depends mostly on the number of generations, which is fixed to $N=200$.

A precise breakeven point for the number of downstream solves depends on the task and the baseline chosen for comparison, and it can range between $500$ solves for the large-scale SUSY task against a robust baseline and around $70000$ for the smaller synthetic tasks against an optimized baseline. One rarely has that many problems sitting and waiting to be solved all at once, but new instances might arrive over the course of days or weeks, and in that case a one-time discovery process running for a few hours (e.g. overnight) can be extremely beneficial.

\section{Discovering Practical LS Algorithms}\label{sec::experiments}

\paragraph{Experimental Setup \& Baselines.} Our discovery framework supports a variety of hardware and programming languages, but in this paper we focus on CPU and Python (see Appendix \ref{app::hardware} for more details). 
We generate synthetic instances $A \in \mathbb{R}^{m \times n}$ and $b \in \mathbb{R}^m$ with $\|A\|_2=\|b\|_2=1$ as follows. First, choose a condition number $\kappa \ge 1$ and fix the optimal residual value to $r^*=0.01$.
\begin{enumerate}
    \item Sample two random matrices with orthonormal columns $U \in \mathbb{R}^{m \times n}$ and $V \in \mathbb{R}^{n \times n}$.
    \item Generate a list of singular values $1=\sigma_1 \ge \sigma_2 \ge ... \ge \sigma_n \ge 0$ decaying exponentially as $\sigma_i = \kappa^{-\frac{i-1}{n-1}}$, and set $\Sigma = Diag(\sigma)$. Set $A = U \Sigma V^\top$.
    \item Let $w \in \mathbb{R}^n$ be a random vector on the unit sphere, $u$ be a unit vector orthogonal to the column space of $A$, and set $b = r^* u + \sqrt{1-(r^*)^2} Uw$. This implies that $\|b\|=1$, the optimal residual is $r^*$, and $x^* = \sqrt{1-(r^*)^2} V \Sigma^{-1} w$.\footnote{This ensures $Ax^*$ is isotropic on the column space of $A$, so iterative algorithms do not have unusually fast convergence due to alignment between solution and singular vectors. As usual, there will be some floating point error in the stored $x^*$ that degrades with $\kappa$.}
\end{enumerate}
We consider several synthetic tasks, varying different parameters such as matrix size, condition number, and target accuracy. We compare the discovered programs with several state-of-the-art baselines (Appendix \ref{app::baselines}):
\begin{itemize}
    \item \textbf{QR}: Householder-QR decomposition based on the LAPACK driver ``\texttt{gels}'', which is the robust full-rank reference solver for the task. 
    \item \textbf{Cholesky}: Direct solver for the normal equations based on the Cholesky decomposition of the Gram matrix, which is the fast direct baseline.
    \item \textbf{Blendenpik} \citep{avron2010blendenpik}: We set the embedding dimension $d=2n$ and warm-start the iterative algorithm to the sketch-and-solve initialization.
    \item \textbf{LSQRsparse}: A sparse sketch followed by preconditioned LSQR. This replaces the subsampled trigonometric sketching operator used in Blendenpik with an OSNAP sketch of dimension $d=3n$ and column sparsity $s=8$ \citep{nelson2013osnap}. This baseline is inspired by recent randomized solvers \citep{epperly2026fast, huh2026adaptive}, but it is more suitable to the practical settings discussed in this paper, see Appendix \ref{app::baselines} for more details. The initial program used for all tasks is very similar to this baseline (Appendix \ref{app::method}).
\end{itemize}
Throughout all our synthetic tasks, we use $20$ instances in the training phase, and we report the average runtime of the discovered algorithm on a disjoint set of $10$ test instances.
As discussed in Section \ref{sec::shinkalgodisc}, we focus on relative residual ratio $f_{\text{res}}$ as our primary notion of accuracy, but in some practical applications one might be interested in forward or backward error. We explore the robustness of our discovered algorithms to other notions of accuracy in Appendix \ref{app::fwdbwd}.

As a real-world application, we also apply our discovery framework to the SUSY dataset \citep{susy_279}, which is used as a benchmark to detect exotic particles from high-energy collisions in particle physics \citep{baldi2014searching}. We consider the kernel regression problem following \citet{epperly2024fast} with $n=2000$ features. This is a naturally dense tall least squares problem, and discovering faster tailored algorithms might help researchers classify more collisions efficiently. \footnote{We omit the error bars because all runtimes are within $1-3\%$ of the average for all synthetic tasks. More fine-grained performance statistics are discussed in Appendix \ref{app::susy_deepdive}.}

\begin{figure}[h]
\begin{center}
    \includegraphics[width=\textwidth]{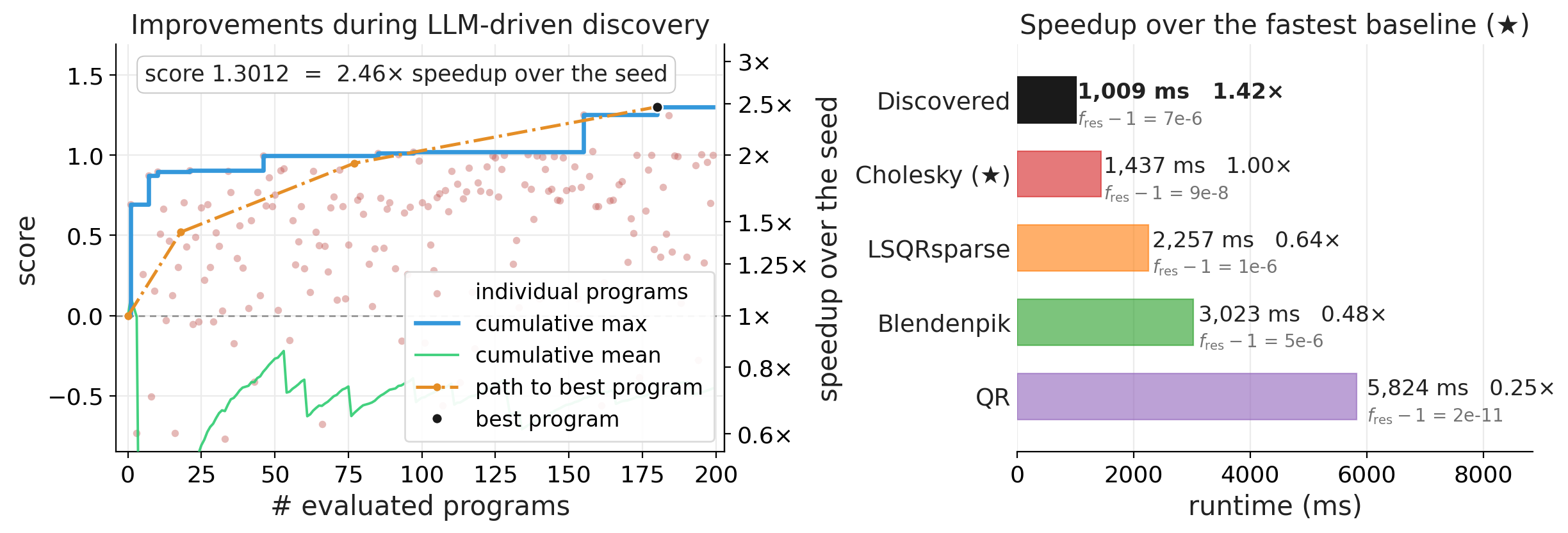}
\end{center}
\caption{\textbf{Left.} Evolution of the search performance scores during the LLM-guided discovery phase. We observe initial rapid runtime performance improvements with sustained minor tuning. \textbf{Right.}  Performance of the discovered algorithm against the baselines. The discovered sketch-and-precondition algorithm outperforms all considered baselines, with $1.42\times$ to $5.77\times$ speedup.}
\label{fig::4.1}
\end{figure}

\paragraph{Single-Task LS Discovery Results.} The first synthetic task consists of matrices with size $m=30000$ and $n=2000$, condition number $\kappa = 10^6$, and target accuracy $\varepsilon=10^{-4}$. This is a relatively large ill-conditioned task which needs to be solved with medium accuracy. During the exploration phase, the LLM-guided framework proposes and evaluates many programs, with fast improvement in the first few iterations but also sustained progress in the late stages (see \Cref{fig::4.1}, left). The discovered algorithm is a sketch-and-precondition iterative algorithm, but differs substantially from the initial one, and it is $1.42\times$ faster than the best baseline on the test set (see Figure \ref{fig::4.1}, right). It uses a sparse sketching operator of dimension $3n$ and $1$ non-zero per column, followed by a Cholesky preconditioner which falls back to QR if the Gram matrix loses positive definiteness, and then preconditioned Conjugate Gradient (CGLS) instead of LSQR, together with many low-level optimizations. Moreover, the stopping rule also differs from the one typical in CGLS (Appendix \ref{app::analysis}).

\begin{figure}[h]
\begin{center}
    \includegraphics[width=\textwidth]{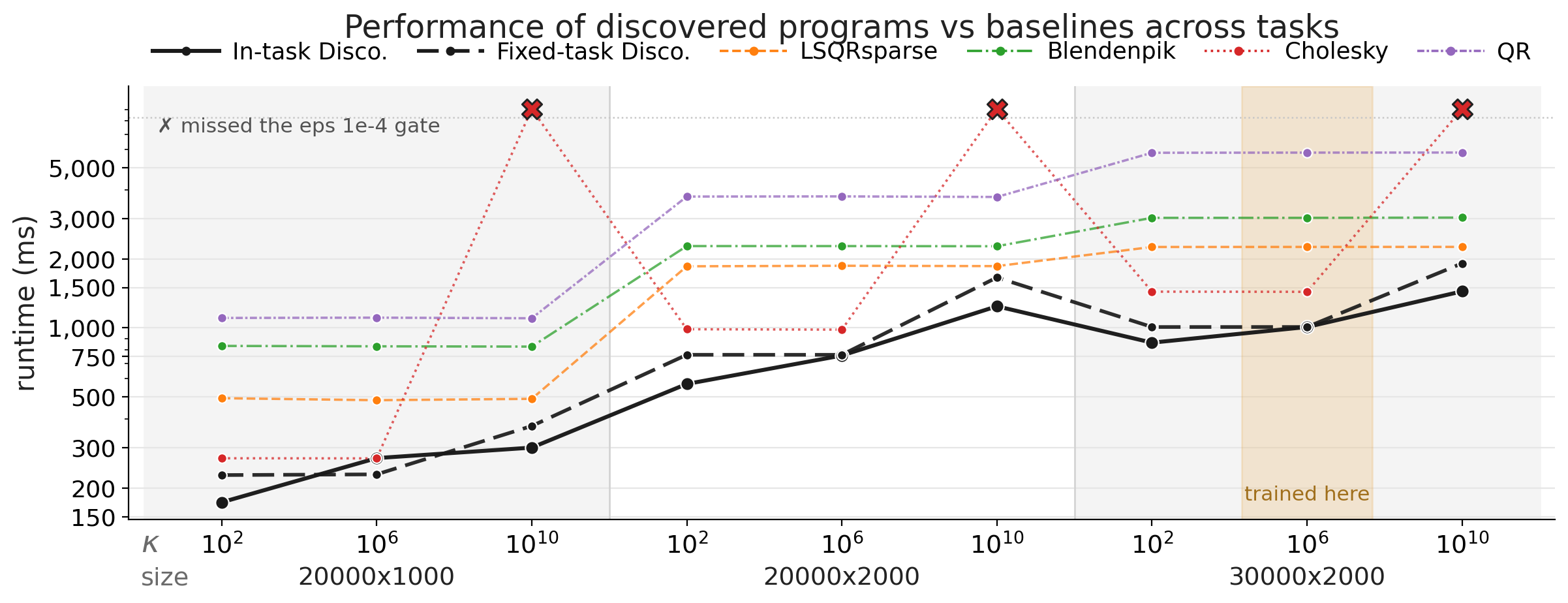}
\end{center}
\caption{Test-set performance comparison between baselines and in-task discovered algorithms (solid black line) and a fixed discovered algorithm (dashed black line), across a variety of tasks. Our framework discovers comparable or faster algorithms across all tasks, and these algorithms generalize well outside the task they were trained on.}
\label{fig::4.2}
\end{figure}

\paragraph{Multi-Task LS Discovery Results.} To study the flexibility of our method, we considered a grid of neighboring tasks, where we vary matrix size and condition number, with $(m,n) \in \{(20000,1000), (20000,2000), (30000,2000)\}$ and $\kappa \in \{10^2,10^6,10^{10}\}$. The best baseline for all tasks considered here with $\kappa < 10^{10}$ is always Cholesky on the normal equations, although this is not a robust solver because it fails to reach the desired accuracy at large condition numbers. The best baseline for ill-conditioned matrices is the randomized preconditioner from a sparse sketch. Our framework discovers faster algorithms for all these tasks, which are up to $1.73\times$ faster than the best baseline for size $(20000, 2000)$ and condition number $10^2$. This is shown in Figure \ref{fig::4.2}, together with the performance of the program discovered at size $(30000,2000)$ and condition number $10^6$. This shows that the discovered algorithms can be empirically robust and perform better than the baselines across different tasks.

As a third synthetic experiment, we fix matrix size $(20000,2000)$ and condition number $10^6$, and vary the target accuracy level $\varepsilon \in \{1, 10^{-2}, 10^{-4}, 10^{-6}, 10^{-8}, 10^{-10}\}$. Again, we observe a transition where direct Cholesky is the fastest baseline when the task solvability is far from the accuracy floor, and then fails as the required accuracy is tighter. The discovered algorithms tie or beat every baseline in every task again, as shown in Figure \ref{fig::1}.

\paragraph{Real-World LS Discovery Results.} Beyond the previously studied synthetic tasks, we investigate the usefulness of our proposed workflow in an applied real-world setting. Using the SUSY dataset \citep{susy_279} for high-energy particle physics, we form kernel regression least squares tasks by selecting $10$ disjoint blocks of $100,000$ rows used for training, which correspond to matrices of size $100,000 \times 2000$. We then evaluate the algorithm on three disjoint sets of $1,000,000$ rows. Since the exact minimizer $x^*$ is not known for these instances, we use the result of a QR direct solve as a reference. On a medium accuracy task $\varepsilon = 10^{-4}$, our method discovers an algorithm which finds a solution $2\times$ faster than the next best baseline, while on a high accuracy task $\varepsilon = 10^{-12}$ the speedup is $1.52\times$, see Table \ref{tab::susy} for more details. 
The test matrices considered here have condition numbers between $9.5\times 10^6$ and $1.06 \times 10^7$, coherence between $46$ and $122$, and the following average row-leverage profile.\footnote{This is a much more skewed leverage profile than the synthetic matrices previously considered.}
\begin{equation*}
    (q0.01, q0.1, q0.5, q0.9, q0.99) = (0.15, 0.25, 0.64, 5.3, 8.8)
\end{equation*}
The baselines tend to over-solve the medium-accuracy task compared to the discovered algorithm, but in Appendix \ref{app::susy_deepdive} we show that the early stopping alone does not explain its advantage.

\begin{table}[h]
\centering
\setlength{\tabcolsep}{6pt}
\caption{Performance of discovered algorithms on instances of $1,000,000$ rows from the SUSY dataset. The error $f_{\text{res}} -1$ is computed relative to the reference solution found by the QR direct solver, hence the $0$ in the corresponding row. The reported error and runtime are the averages over the test set, and Cholesky fails to reach the target accuracy for $\varepsilon=10^{-12}$ (red).}
\label{tab::susy}
\begin{tabular}{l rrr rrr}
\toprule
 & \multicolumn{3}{c}{$\varepsilon = 10^{-4}$} & \multicolumn{3}{c}{$\varepsilon = 10^{-12}$} \\
\cmidrule(lr){2-4} \cmidrule(lr){5-7}
Program & time (s) & speedup & error & time (s) & speedup & error \\
\midrule
Discovered & \textbf{20.7} & $2.00\times$ & $2.4\times10^{-5}$ & \textbf{48.3} & $1.52\times$ & $1.1\times10^{-15}$ \\
LSQRsparse \textbf{($\star$)} & 41.4 & $1.00\times$ & $9.4\times10^{-7}$ & 73.7 & $1.00\times$ & $4.6\times10^{-14}$ \\
Blendenpik & 94.0 & $0.44\times$ & $3.4\times10^{-6}$ & 145.4 & $0.51\times$ & $1.1\times10^{-13}$ \\
Cholesky & 45.6 & $0.91\times$ & $8.2\times10^{-11}$ & \textcolor[HTML]{C44E52}{45.4} & \textcolor[HTML]{C44E52}{--} & \textcolor[HTML]{C44E52}{$8.2\times10^{-11}$} \\
QR & 242.7 & $0.17\times$ & $0^{*}$ & 243.2 & $0.30\times$ & $0^{*}$ \\
\bottomrule
\end{tabular}
\end{table}

\section{Analysing the Discovered LS Algorithms}\label{sec::analysis}

Next, we provide an in-depth study of the algorithms discovered by our methodology, highlighting the sources of efficiency gains. Moreover, we compare our LLM-guided discovery framework with hyperparameter tuning, showing that LLM search can outperform hyperparameter tuning, but the two can be combined to find even faster algorithms. 

\begin{figure}[h]
\begin{center}
    \includegraphics[width=\textwidth]{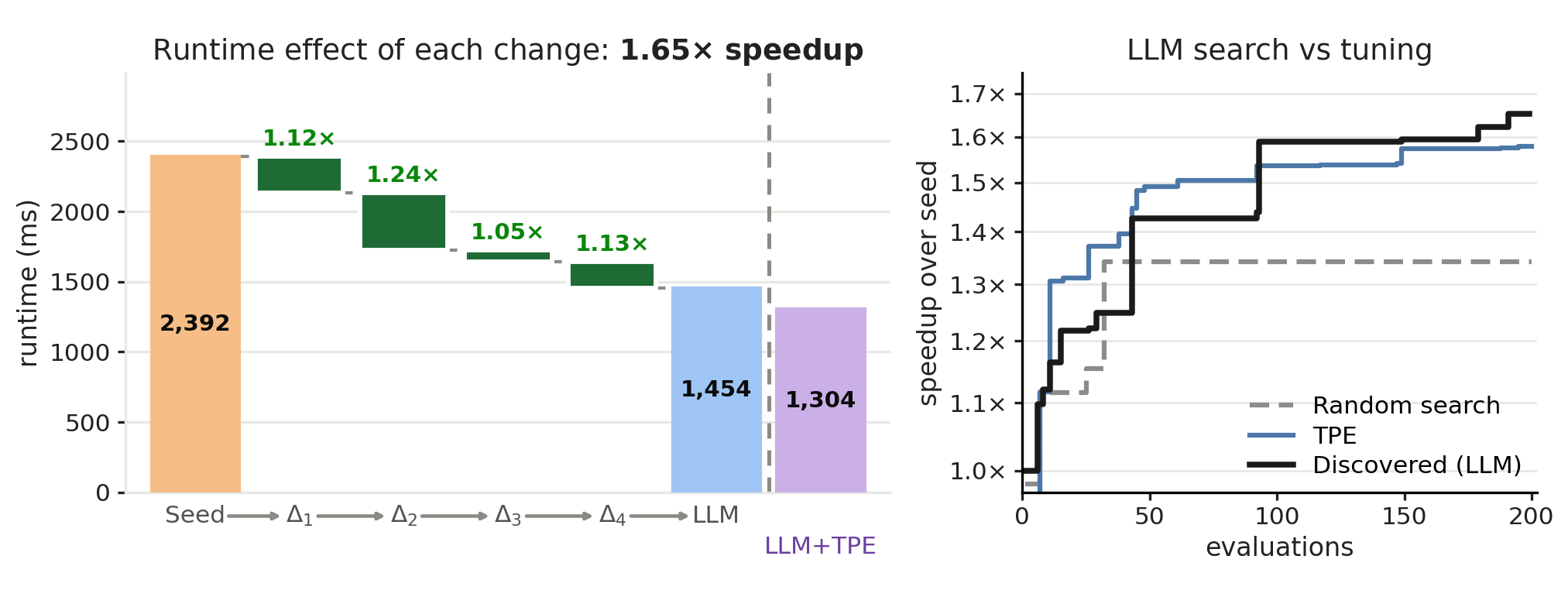}
\end{center}
\vspace{-1em}
\caption{\textbf{Left.} Breakdown of the runtime improvement achieved by the discovered algorithm. $\Delta_1$ changed the termination condition \textit{atol} from $\varepsilon$ to $0.16\sqrt{\varepsilon}$, $\Delta_2$ changed the sketching family from OSNAP with 8 nonzeroes per column to row-sampling with row-norms probabilities, $\Delta_3$ lowered the sketch dimension from $3n$ to $1.9n$, and $\Delta_4$ includes some code optimizations. Performing hyperparameter tuning via TPE on the discovered algorithm can further improve its performance by $1.11\times$. \textbf{Right.} Comparison between LLM-guided discovery (black) and hyperparameter optimization starting from the same initial program for this task. 
We considered two hyperparameter search methods: $200$ iterations of TPE (blue) or  taking the best of $200$ random valid configurations (grey).
}
\label{fig::breakdown_vs_hyper}
\end{figure}

\paragraph{Algorithm Analysis for a Fixed Task.} We consider the algorithm discovered in the task with matrices of size $(30000,2000)$, condition number $10^{10}$, and accuracy target $10^{-4}$. This remains a randomized preconditioner, but most pieces of the initial program were changed, including the termination condition on the backward error ($\Delta_1$), the sketch family ($\Delta_2$), the sketch dimension ($\Delta_3$), together with some low-level code optimization ($\Delta_4$). See Figure \ref{fig::breakdown_vs_hyper} for a detailed breakdown of each contribution, while we defer the full description of each program to Appendix \ref{app::analysis}, which also contains the analysis of algorithms discovered for different tasks. One particularly interesting discovery is a randomized preconditioner with a new stopping criterion which is mathematically sound and provides a speedup over the traditional backward error rule (see Appendix \ref{app::analysis}, Eq. \ref{eq::newbound2}).

\paragraph{Program Evolution versus Hyperparameter Optimization.}
Comparing arbitrary program evolution with hyperparameter optimization is generally delicate because many types of mutations are difficult to represent as hyperparameters. However, many randomized iterative algorithms for least squares naturally lend themselves to it \citep{cho2025surrogate}, so we compared our framework with random hyperparameter search and more advanced hyperparameter tuning using the tree-structured Parzen estimator \citep[TPE, ][]{bergstra2011algorithms} implemented via Optuna \citep{akiba2019optuna} on the task of the previous paragraph. Our framework outperformed both random search and TPE, and it kept producing improvements at later stages of the exploration. We also considered combining the two approaches, where we run a few iterations of TPE starting from the program discovered by our framework, and this led to even faster algorithms; see Figure \ref{fig::breakdown_vs_hyper}. This suggests that LLM-guided evolution can make larger ``jumps'' in late iterations, which helps prevent search stagnation typical of Bayesian hyperparameter tuning, but the local optimization provided by hyperparameter tuning can further refine these solutions. However, not all discovered algorithms might naturally lend themselves to hyperparameter tuning, for example direct algorithms do not expose a comparable set of tunable hyperparameters as randomized iterative methods. Another drawback of hyperparameter tuning is that the search space is constrained by a user-specified family of algorithms defined by the hyperparameters chosen, while LLM-guided search is unrestricted and much more flexible. The hyperparameter configurations considered in the search is much broader than in previous works, and was constructed using extensive expert knowledge (Appendix \ref{app::hyper}) to make tuners as competitive as possible with our framework, which would otherwise trivially outperform tuners by a large margin on these tasks.

\paragraph{Program Evolution versus Random Sampling.} We compare our framework based on LLM-guided evolutionary search and a \textit{best-of-N} LLM sampling strategy. Both methods achieve a significant improvement over the initial program, but the evolutionary approach used by our framework yields a $1.16\times$ speedup over random sampling from LLMs, thanks to the deeper exploitation of previously proposed ideas (Figure \ref{fig::4.3}, right). Detailed results are provided in Appendix \ref{app::shinka_vs_restarts}.
\paragraph{Reproducibility and further comparison with autotuners} So far, every discovered program is the result of a single run of our framework described in Section \ref{sec::shinkalgodisc}. However, LLMs are inherently stochastic and can produce very different answers to the same queries. Independent runs of our framework on the same task consistently produced significant runtime improvements over the strongest baselines. The exact discovered algorithms are slightly different from run to run, but their performance is similar, which is desirable in practice. More details on reproducibility of our discovery framework are provided in Appendix \ref{app::reproducibility}. We also investigate the reproducibility of hyperparameter tuning. The final performance gains of LLM-guided evolution over a well-specified tuning are modest yet significant, however hyperparameter tuning suffers from mis-specification much more, and requires more expert knowledge to produce successful results. 
In Appendix \ref{app::reproducibility} we also provide a comparison with Gaussian-Process surrogate-based autotuning, where the performance gains of our framework are even greater. This is possibly due to the difficulty of Bayesian optimization in handling high-dimensional search spaces with categorical variables \citep{eggensperger2013towards}, such as the one of Table \ref{tab::hyper}. 
\paragraph{Out-of-Distribution Generalization.} To understand the out-of-distribution generalization capabilities of the algorithms discovered by our framework, we consider a task with a mix of two families of instances evenly split between $((m,n),\kappa) \in \{((20000,1000), 10^{10}), ((30000,2000), 10^{2})\}$. We evaluate the discovered algorithm over a grid of instances with different size and condition number, and different target accuracy levels (Figure \ref{fig::4.3}, left). Runtime speedups transfer strongly across size and conditioning, where the algorithm outperforms all baselines, but the target accuracy generalizes only on the looser side, with the algorithm failing to reach the target accuracy at $\varepsilon = 10^{-6}$.  We also studied generalization capabilities to tasks which have larger structural differences from the distribution of instances used for discovery, including spectral distribution, coherence levels, and optimal residual scales (Appendix \ref{app::OOD}). In particular, we show that the discovery can be fragile to severe task mis-specification, but it can be largely repaired when the mis-specification is milder.

\begin{figure}[h]
\begin{center}
    \includegraphics[width=0.8\textwidth]{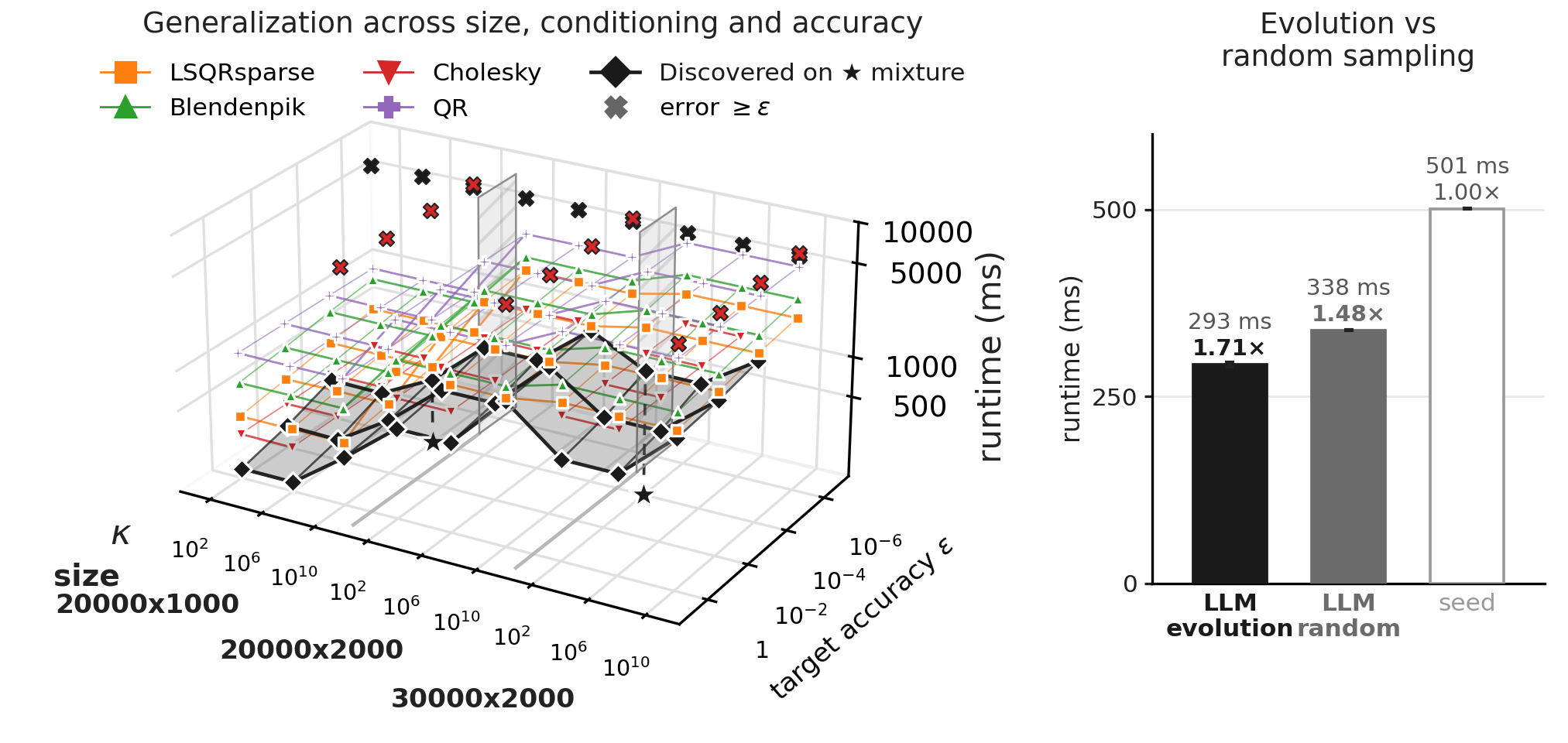}
\end{center}
\vspace{-1em}
\caption{\textbf{Left.} Generalization of a discovered LS algorithm across a grid of different size, conditioning, and accuracy targets. The algorithm is trained at accuracy level $10^{-4}$ and a mixture of instances of size and conditioning moving in opposite directions, and its generalization is one-sided, working well at lower accuracy but failing at higher accuracy. The left axis flattens the size and condition number axes, making it possible to visualize this 4-dimensional data in 3D. \textbf{Right.} Comparison between our framework (LLM evolution) and best-of-N sampling (LLM random) on a task with size $(20000,1000)$, $\kappa=10^{10}$, $\varepsilon=10^{-4}$, yielding $1.71\times$ and $1.48\times$ speedup over the seed.}
\label{fig::4.3}
\end{figure}


\section{Conclusion}\label{sec::conclusions}

\paragraph{Summary.} We introduced a framework for using LLM-guided search to discover faster algorithms tailored to practical least squares applications. Our results demonstrate significant runtime speedups over robust baselines and default implementations with minimal task-specific input, offering a way to automate the tedious hand-tuning required to speed up stable algorithms for LS tasks.
Our framework for evaluating practical-case algorithmic tasks is flexible and can be easily modified beyond the least squares application discussed in this paper.

\vspace{-0.2em}

\paragraph{Limitations.} Our framework requires pre-computing solutions to some instances and repeatedly solving many least squares problems during evaluation, so the benefits of the faster discovered algorithms become tangible in applications which require several downstream or repeated solves. The algorithms discovered are often a recombination of LS techniques known to experts, whose runtime improvement can partly be achieved by well-specified hyperparameter tuning on a large space.
Finding a genuinely novel high-level algorithmic discovery remains challenging. The algorithms discovered lack provable guarantees beyond their success on task instances, so they should be carefully inspected.

\vspace{-0.2em}

\paragraph{Future Directions.} Given the flexibility of our framework, it would be very interesting to apply it to other tasks beyond dense least squares, for example to sparse settings or to other linear algebraic problems. Further understanding and improving the LLM-guided exploration framework to discover new provably stable algorithms is also an interesting research direction. Moreover, our framework currently optimizes a scalar runtime objective subject to a scalar feasibility constraint, and extending it to optimize a Pareto front of runtime and different accuracy metrics is another interesting future direction.

\vspace{-0.2em}

\paragraph{Broader Impact \& Ethical Considerations.} We hope our work can be useful to practitioners who need to quickly solve many large-scale least squares problems. The algorithms discovered for the tasks considered in this paper and our code are publicly available in the following repository: \href{https://github.com/SakanaAI/AutoDiscoLS}{https://github.com/SakanaAI/AutoDiscoLS}.

\subsubsection*{Acknowledgments}

We thank Edoardo Cetin, Antares Chen, Nick Lindsay and Daniel Spielman for helpful discussions on different parts of this project. This work was supported, in part, by NSF CCF-2112665 (TILOS) and by ONR Award N00014-24-12611.

\bibliography{references}
\bibliographystyle{iclr2027_conference}

\newpage

\appendix

\startcontents[appendix]
\section*{Appendix table of contents}
\printcontents[appendix]{}{1}{\setcounter{tocdepth}{2}}

\section{Hardware and Software Packages}\label{app::hardware}

All experiments were performed on a shared node inside a SLURM-managed partition with 2TB of DDR5 memory and Intel Xeon Platinum 8481C 2.7GHz CPU with 104 cores. Up to 60 cores were used during the LLM-guided search phase to parallelize candidate evaluation, but all timings on test instances presented in the paper were measured with exclusive use of the node and 32GB of memory. We used CPython 3.14.5 with NumPy 2.4.6, SciPy 1.17.1, OpenBLAS 0.3.30, and ShinkaEvolve 0.0.7. 

\section{More details on the methodology}\label{app::method} 

\subsection{Evaluation function}
The log geometric-mean speedup function (\ref{eq::shinka-score-corr}) is designed to reduce runtime noise (hence the median) but also preserve signals coming from variations across instances (hence the mean).

Let $S \subseteq I_{\text{train}}$, and define the restriction of (\ref{eq::shinka-score-corr}) on $S$ as
\begin{equation*}
    f(\mathcal{A};S) = -\frac{1}{|S|}\sum_{i\in S} \log_2 \left( \text{median}_{r=1,...,R} (\rho(i,r)) \right)
\end{equation*}
and let $\sigma(x) = \frac{1}{1+e^{-x}}$ be the sigmoid function. Define the infeasible scores for (\ref{eq::shinka-score}) as
\begin{equation*}
    g(\mathcal{A},S) = \begin{cases}
    K - |S| - \sigma(f(\mathcal{A};S)) & \text{if } |S| > 0 \\
    K - |S|  & \text{if } |S| = 0 
    \end{cases}
\end{equation*}
This objective function was chosen so that programs with more correct instances are always better than programs with fewer ones, and speed differentiates programs that solve the same number of instances. The floor $\tau$ is added to separate the feasible and infeasible region. Throughout all our experiments, we fix $\tau=-5$ and $R=3$.

The function (\ref{eq::shinka-score}) is designed so that a feasible initial program $\mathcal{A}_0$ always has score $0$, with some tiny fluctuations due to runtime measurement noise, while programs with score greater than $0$ are always feasible and faster. The score function is chosen so that $\text{score}(\mathcal{A})=1$ means that the algorithm $\mathcal{A}$ takes half the time as $\mathcal{A}_0$, $\text{score}(\mathcal{A})=-1$ means double the time, and in general $\text{score}(\mathcal{A})=s$ means a $2^{-s}$ fraction of the time.

To reduce runtime noise, we adopt the following measures.
\begin{itemize}
    \item We normalize the state of the machine/caches by running a fixed program before any timed evaluation.
    \item We interleave the initial program and the current program evaluations, i.e. we compute $t_{\mathcal{A}}(1,1),t_{\mathcal{A}_0}(1,1), t_{\mathcal{A}}(1,2),t_{\mathcal{A}_0}(1,2), ..., t_{\mathcal{A}}(K,3),t_{\mathcal{A}_0}(K,3)$ in this order.
    \item We pin processes used in the evaluation to physical cores.
\end{itemize}

\subsection{ShinkaEvolve parameters}
All the experiments in this paper used the same initial program (\Cref{code::initial}) and the same configuration of ShinkaEvolve parameters (\Cref{tab::shinka}). All runs completed $200$ generations using around $\$10$ per task despite the larger cost budget, but one can lower the cost even further by using a cheaper LLM ensemble. This might lead to lower sample efficiency to reach comparable algorithms, but it can still produce improvements over the baselines for a fraction of the cost.

\begin{lstlisting}[
  language=Python,
  caption={Task \texttt{msg} parameter used as \texttt{task\_sys\_msg} inside ShinkaEvolve. The only part that changes across different tasks is the one on the target accuracy, where we report the desired value and label $\varepsilon \in \{1,0.01\}$ as ``low-precision'', $\varepsilon \in \{10^{-4}, 10^{-6}\}$ as ``moderate-precision'', and $\varepsilon \in \{10^{-8}, 10^{-10}\}$ as ``high-precision''.},
  label={code::msg},
  % basicstyle=\ttfamily\tiny,  % Change font size here
]
"""
You are an expert in numerical linear algebra and randomized algorithms.

TASK: write ONE Python function that solves the overdetermined least-squares problem  min_x ||A x - b||_2
(A is a real numpy array of shape (m, n) with m >> n; b is a real numpy vector of shape (m,)) as FAST as possible in measured wall-clock time. 

`eps` is the target residual accuracy ||A x - b|| <= (1 + eps) * min_x ||A x - b||. 
For this task, we are interested in moderate-precision solutions, so eps is on the order of 1e-4.  However, it is a RUN-TIME PARAMETER, not a constant, and you cannot check that condition without solving the problem, so DECIDING WHEN TO STOP IS PART OF THE PROBLEM.

Build a stopping test out of what you can actually see -- A, b, your current iterate -- and make it cheap, because whatever it costs is timed, but also accurate enough (failing to solve the problem to the desired accuracy results in a large penalty).

You receive input data in double-precision floats, but feel free to try lower precision arithmetic to make the algorithm faster.

GOAL: minimize the wall-clock time to a solution that reaches the accuracy on EVERY problem. Be fast AND dependable -- a program that fails an instance, or that is occasionally far slower than usual, is not good enough.

Feel free to use optimized linear algebra libraries like numpy and scipy. 
Try sketching or a combination of randomized and direct methods in any way you see fit.

Because the objective is REAL wall-clock time, cheap-but-good beats theoretically-optimal: prefer whole-array numpy BLAS-3 operations when possible, avoid recomputing, and minimize total time to an accepted x.
"""
\end{lstlisting}

\begin{lstlisting}[
  language=Python,
  caption={Initial program used by our discovery framework},
  label={code::initial},
  basicstyle=\ttfamily\tiny,  % Change font size here
]
# EVOLVE-BLOCK-START
import numpy as np
from scipy import linalg as sla
from scipy.sparse import csr_matrix
from scipy.sparse.linalg import LinearOperator, lsqr

def _osnap(m, d, s):
    cols = np.repeat(np.arange(m), s)
    rows = np.random.randint(0, d, size=m * s)
    signs = np.random.choice((-1.0, 1.0), size=m * s) / np.sqrt(s)
    return csr_matrix((signs, (rows, cols)), shape=(d, m))

def solve(A, b, eps, iter_budget):
    m, n = A.shape

    # --- Preconditioner: sketch A (OSNAP, d=3n rows), thin-QR of the sketch, M = R^{-1} (right-preconditioner).
    d = 3 * n
    S = _osnap(m, d, 8)
    SA = S @ A
    Sb = S @ b
    R = np.linalg.qr(SA, mode="r")
    M = sla.solve_triangular(R, np.eye(n))

    # --- Warm start: sketch-and-solve, so LSQR only iterates on the residual.
    x0 = M @ (M.T @ (SA.T @ Sb))
    b_eff = b - A @ x0
    if np.linalg.norm(b_eff) <= 1e-20:
        return x0

    # --- Preconditioned LSQR on A @ M, stopping on its own backward-error test at tolerance eps.
    op = LinearOperator((m, n), matvec=lambda v: A @ (M @ v),
                        rmatvec=lambda u: M.T @ (A.T @ u), dtype=np.float64)
    y = lsqr(op, b_eff, atol=eps, btol=eps, conlim=0, iter_lim=iter_budget)[0]
    return x0 + M @ y

# EVOLVE-BLOCK-END

# Fixed harness (timing + accuracy check), hidden from the search: it cannot bypass the timing or the grading.

from program_run import run_solve as _run

def run_solve(A, b, nnz, eps, iter_budget, run_idx, r_star, seed):
    return _run(solve, A, b, nnz, eps, iter_budget, run_idx, r_star, seed)
\end{lstlisting}

\begin{table}[h]
    \centering
    \caption{ShinkaEvolve parameters that differ from the default ones. The number of generations is 201 because it includes the seed program.}
    \label{tab::shinka}
    \begin{tabular}{@{}ll@{}}
    \toprule
    Parameter & Value \\
    \midrule
    \multicolumn{2}{@{}l}{\texttt{EvolutionConfig}} \\
    \addlinespace[2pt]
    \texttt{num\_generations} & 201 \\
    \texttt{max\_api\_costs} & 20 \\
    \texttt{llm\_models} & \texttt{gpt-5.3-codex} \\
     & \texttt{gpt-5.6-terra} \\
     & \texttt{gpt-5.6-sol} \\
     & \texttt{gemini-3.5-flash} \\
     & \texttt{gemini-3.1-pro-preview} \\
    \texttt{llm\_kwargs.reasoning\_efforts} & low, medium, high \\
    \texttt{llm\_kwargs.max\_tokens} & 64000 \\
    \texttt{novelty\_llm\_models} & \texttt{gpt-5.6-luna} \\
    \texttt{code\_embed\_sim\_threshold} & 0.92 \\
    \midrule
    \multicolumn{2}{@{}l}{\texttt{DatabaseConfig}} \\
    \addlinespace[2pt]
    \texttt{num\_islands} & 4 \\
    \texttt{enable\_dynamic\_islands} & True \\
    \texttt{stagnation\_threshold} & 50 \\
    \texttt{migration\_rate} & 0.1 \\
    \texttt{migration\_interval} & 40 \\
    \texttt{archive\_size} & 80 \\
    \texttt{archive\_selection\_strategy} & crowding \\
    \texttt{num\_archive\_inspirations} & 4 \\
    \texttt{num\_top\_k\_inspirations} & 2 \\
    \texttt{parent\_selection\_lambda} & 3 \\
    \midrule
    \multicolumn{2}{@{}l}{\texttt{JobConfig}} \\
    \addlinespace[2pt]
    \texttt{numeric\_threads\_per\_job} & 1 \\
    \midrule
    \multicolumn{2}{@{}l}{\texttt{ShinkaEvolveRunner}} \\
    \addlinespace[2pt]
    \texttt{max\_proposal\_jobs} & 12 \\
    \texttt{max\_evaluation\_jobs} & 6 \\
    \texttt{max\_db\_workers} & 1 \\
    \bottomrule
    \end{tabular}
\end{table}

\section{Tradeoff between LLM-guided evolution and random sampling}\label{app::shinka_vs_restarts}

In this section, we compare the approach described in Section \ref{sec::shinkalgodisc} with random sampling from an ensemble of LLMs. We also consider hybrid approaches where we randomly restart $k$ ShinkaEvolve instances of $N / k$ generations each. The extreme cases of $k=1$ and $k=N$ recover ShinkaEvolve and random LLM sampling, respectively. For settings where we use multiple restarts, we set aside a validation instance set $I_\text{valid}$ different from $I_{\text{train}}$ and $I_{\text{test}}$, which is used to co-time the discovered algorithms from the independent runs and select the best one. This is non-trivial because the score during the evolutionary phase is noisier than during validation and testing due to contention coming from parallel evaluations, so it should not be used to select the best algorithms. In fact, we observed a different validation ranking compared to the training one on multiple occasions.

We fix $N=200$ throughout and consider $k=1,2,5,10$, and also consider two versions of random sampling. The first uses \texttt{msg} as in \Cref{code::msg} and provides the initial program (\Cref{code::initial}) to the LLM prompt, together with the same public metrics as the evolutionary runs (i.e. this is equivalent to $k=200$). The second one is labeled ``naive'' because it provides the most basic seed algorithm (the default \texttt{scipy.linalg.lstsq} or direct Cholesky) and more limited context. This was meant to study the effect of prompting on the LLMs' ability to discover fast least squares algorithms.

\begin{figure}[h]
\begin{center}
    \includegraphics[width=\textwidth]{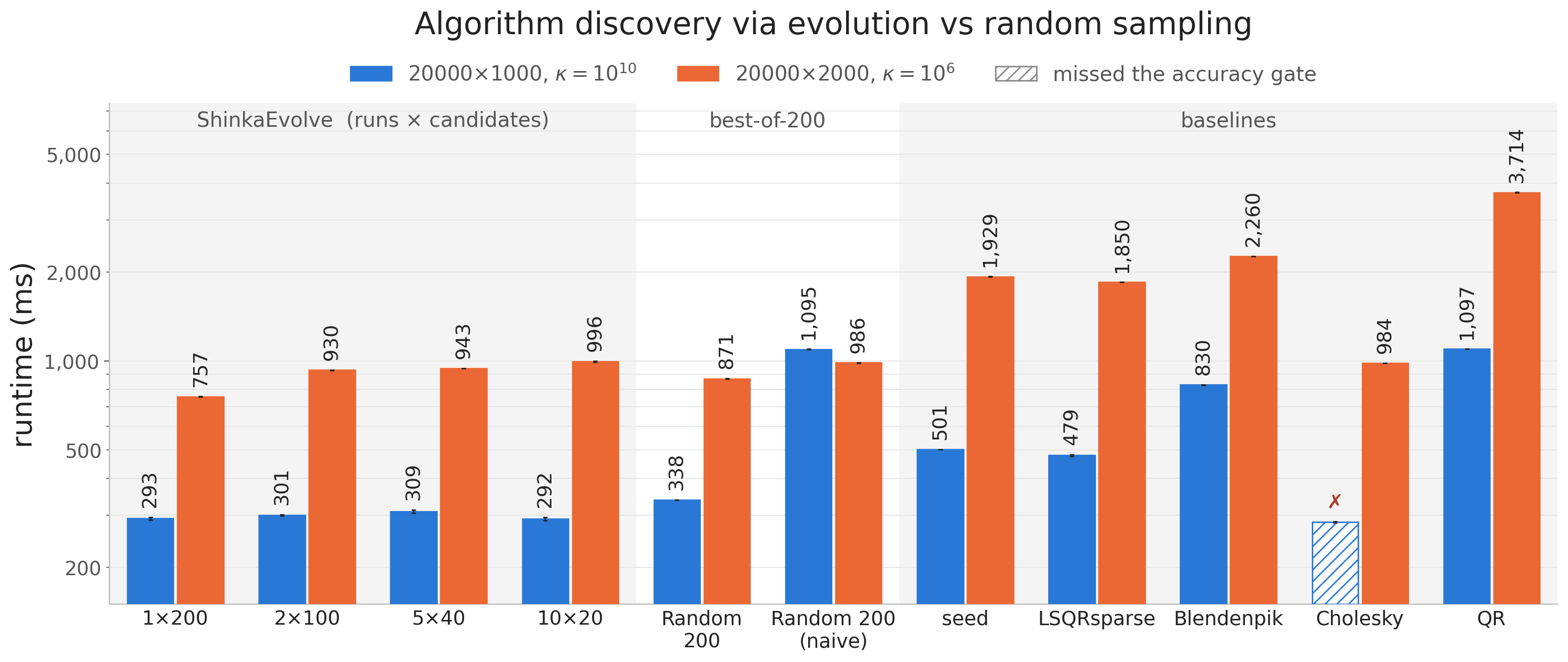}
\end{center}
\caption{Comparison of performance of algorithms discovered under different settings of program refinement (evolution) and random sampling (restarts). On the tasks considered here, the setting with a single LLM-guided evolutionary run outperformed the other ones, but not by a large margin relative to how much most regimes improve over the baselines. On the regimes with best-of-$200$, prompting was essential in exploring beyond the initial program. Direct Cholesky failed to reach the target accuracy for the task with condition number $10^{10}$.}
\label{fig::shinka_vs_restarts}
\end{figure}

We ran all discovery regimes on two tasks, one involving large ill-conditioned matrices (size $(20000,2000)$ and $\kappa=10^6$), and one involving smaller but even more ill-conditioned matrices (size $(20000,1000)$ and $\kappa=10^{10}$), both with target accuracy $\varepsilon=10^{-4}$. Unsurprisingly, prompting has a strong influence on the exploration abilities of LLMs, as evidenced by the ``naive'' runs not leaving the initial program on both tasks. On the other hand, all other discovery regimes exhibited substantial improvements over the initial program and the other baselines. For the very ill-conditioned task, all evolutionary regimes performed similarly, with $k=1$ and $k=10$ as the best ones. The best regime on the other task was $k=1$, outperforming $k=200$ by $1.15\times$ and $k=2$ by $1.23\times$ (\Cref{fig::shinka_vs_restarts}). In light of this, we chose the settings with minimal restarts as our default methodology. This also makes the framework simpler and avoids the need for a validation phase to select the best program among the independent restarts. However, one could have chosen other settings described here and would have likely reached similar results, aside from the naive one. Finding more interesting ways to balance exploration with exploitation for optimizing functions over the space of algorithms remains a very interesting and open research question.

\section{Baselines}\label{app::baselines}

\paragraph{Direct method I: Householder QR} The first baseline is the classical stable solver coming from computing $A=QR$ and solving a triangular system $Rx=Q^\top b$, where $Q^\top b$ is computed via Householder reflectors instead of storing $Q$ explicitly. We use the LAPACK driver \textit{dgels} and manually implement this instead of relying on the default \texttt{scipy.linalg.lstsq} because that uses an SVD decomposition and it is unnecessarily slow for full-rank least squares problems.

\paragraph{Direct method II: Cholesky on the Gram matrix} This is another classical algorithm that forms the normal equations $A^\top Ax = A^\top b$, which is a positive definite linear system when $A$ has full column-rank, so it can be solved by computing the Cholesky decomposition of $G=A^\top A$, i.e. let $G = L L^\top$ and solve $LL^\top x = A^\top b$ using two triangular solves. This is not as numerically stable as QR, since it computes the Cholesky factor of a matrix with squared condition number $\kappa(G)=\kappa(A)^2$, so it will not solve as many problems as QR to high-precision. However, it is much faster when it is applicable.

\paragraph{Randomized preconditioner I: Blendenpik} This algorithm was invented by \cite{avron2010blendenpik} and it follows the sketch-and-precondition template from Section \ref{sec::background}. The sketching operator $S \in \mathbb{R}^{d \times m}$ is a subsampled randomized Hadamard transform (SRHT)
\begin{equation*}
    S = \sqrt{\frac{m}{d}} RHD
\end{equation*}
where $D \in \mathbb{R}^{m \times m}$ is a diagonal matrix with entries chosen uniformly in $\{-1,+1\}$, $H\in \mathbb{R}^{m \times m}$ is the scaled Walsh-Hadamard matrix, and $R \in \mathbb{R}^{d \times m}$ is a matrix that restricts an $m$-dimensional vector to $d$ randomly chosen (but fixed) directions (see \cite{tropp2011improved} for more details). The sketch dimension is fixed to $d=2n$, which makes the algorithm faster than the default recommendation of $d=4n$ in \cite{avron2010blendenpik}.

\paragraph{Randomized preconditioner II: sparse sketching} This uses the sparse-sign sketching operator of \cite{nelson2013osnap}, where the matrix $S \in \mathbb{R}^{d \times m}$ has $d=3n$ rows and exactly $8$ nonzeroes per column with values chosen independently in $\left\{-\frac{1}{2\sqrt{2}}, +\frac{1}{2\sqrt{2}}\right\}$. The rows are split in $8$ disjoint blocks and each nonzero location is chosen at random inside each block. This is a sparse matrix that can be applied to $A$ much more quickly than the dense sketching operator of Blendenpik.

\begin{figure}[!t]
\begin{center}
    \includegraphics[width=0.74\textwidth]{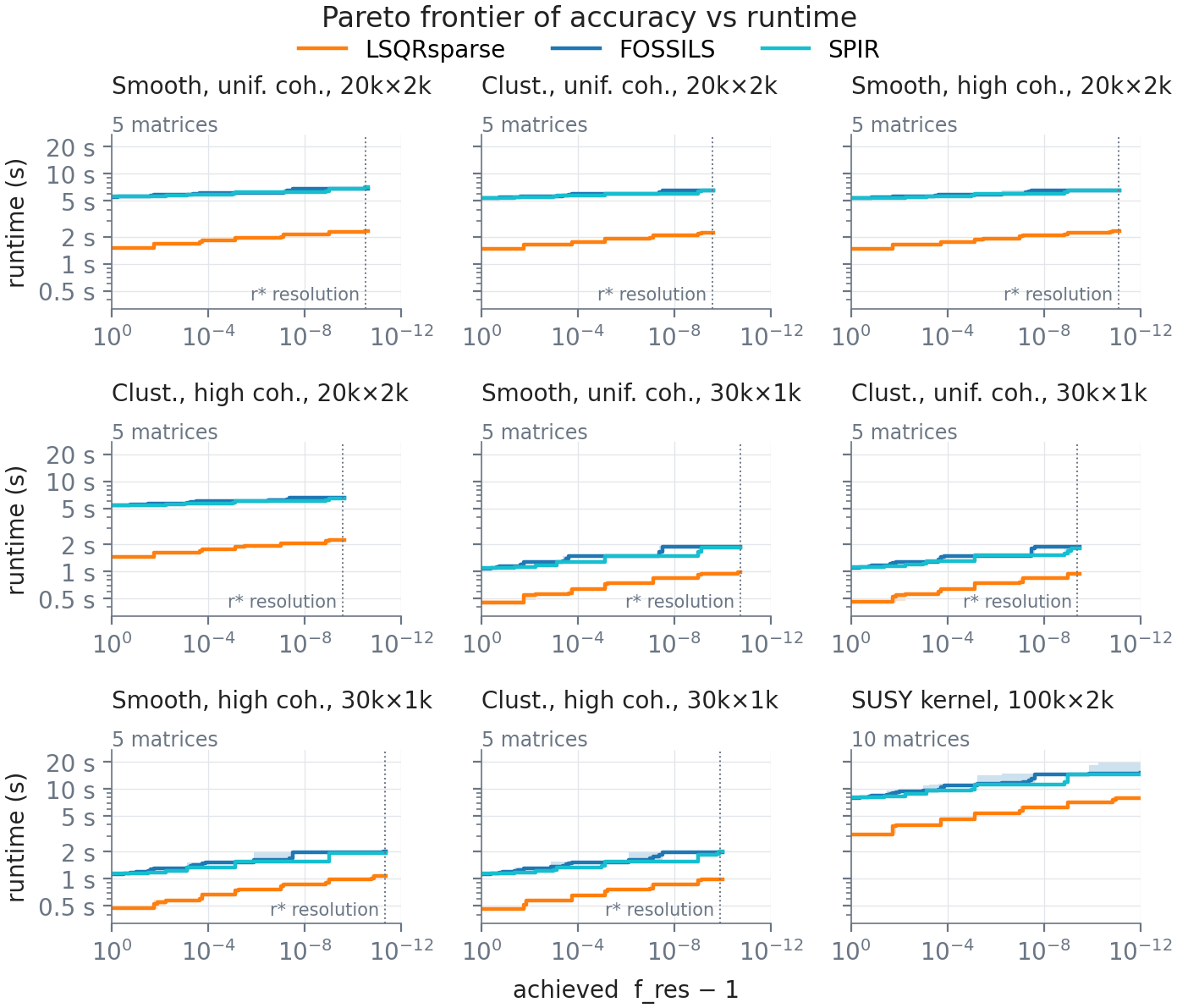}
\end{center}
\vspace{-1em}
\caption{Pareto frontier of runtime and relative residual accuracy for three algorithms on a combination of synthetic datasets using smoothly decaying or clustered (three clusters) singular values, uniform or high-coherence rows, and sizes $30000\times1000$ or $20000 \times 2000$, together with some matrices from the SUSY dataset. We plot the mean and min/max runtime per solver, but this is only noticeable in a few tasks.
}
\label{fig::pareto_baselines}
\end{figure}

\begin{figure}[t]
\begin{center}
    \includegraphics[width=0.74\textwidth]{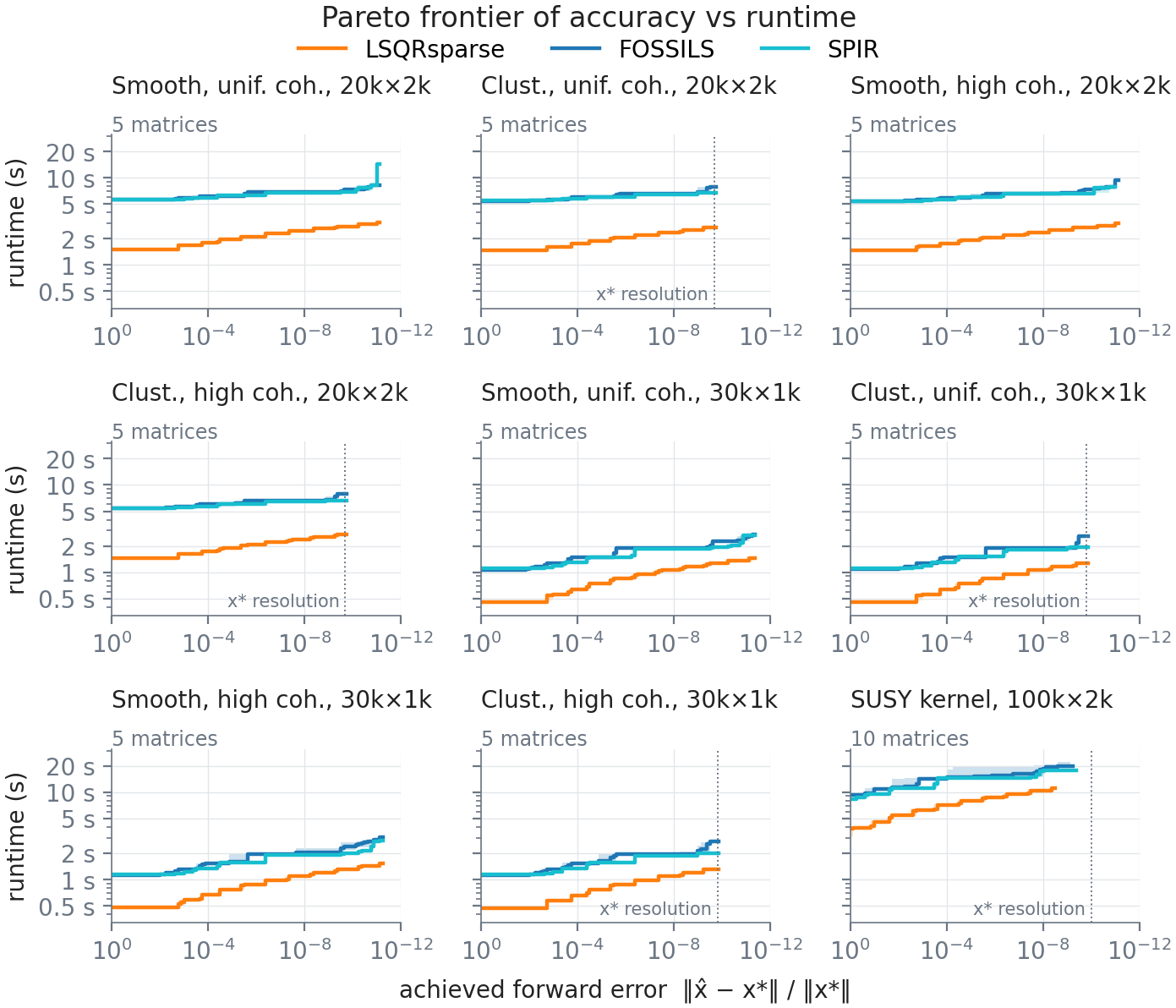}
\end{center}
\vspace{-1em}
\caption{Pareto frontier of runtime and forward error on the same tasks as Figure \ref{fig::pareto_baselines}.
}
\label{fig::pareto_baselines_fwd}
\end{figure}

\paragraph{Termination condition for iterative algorithms} Both randomized preconditioners use the iterative algorithm LSQR on the preconditioned system $\hat{A}=AR^{-1}$, with $R$ coming from the QR factorization of $SA$. This iterative algorithm stops when the scaled normal equation residual $\gamma$ is below a threshold \texttt{atol}\footnote{We set the other termination parameter \texttt{btol} to $0$.}. Let $\hat{r} = \hat{A}x - b$ be the residual at $x$, then this termination condition is
\begin{equation*}
    \gamma = \frac{\|\hat{A}^\top \hat{r}\|}{\|\hat{A}\|\|\hat{r}\|} \le \texttt{atol}
\end{equation*}
This is an upper bound on the backward error in $\hat{A}$ at $x$, and using default parameters might not necessarily reflect our target accuracy in the relative residual ratio $f_{\text{res}}$. Therefore, we tune the termination parameter to 
\begin{equation*}
    \texttt{atol} = 0.03\frac{\sqrt{\varepsilon(2+\varepsilon)}}{1+\varepsilon}
\end{equation*}
This termination condition gives a sufficient condition for $f_{\text{res}} \le 1+ \varepsilon$ in exact arithmetic whenever the preconditioned matrix has condition number less than $30$, which is the case for all applications in this paper.\footnote{In finite precision, there might still be some error introduced by LSQR's running estimate of $\|\hat{A}\|$ through its Frobenius norm.} The term depending on $\varepsilon$ is the result of converting a bound on $\gamma$ to a bound on $f_{\text{res}}$, while the constant $0.03$ is fixed to control the bound on the preconditioned condition number, so lowering it extends the guarantee to even more poorly conditioned matrices.
\paragraph{Comparison with other modern randomized preconditioners} 
The baseline \textbf{LSQRsparse} is inspired by the recommendations in the recent papers \citep{epperly2026fast, huh2026adaptive}, but kept much simpler to reduce robustness and increase efficiency, making it a more practical alternative. This is a non-standard baseline designed using our own expertise from working with these algorithms in practice. To motivate this choice, we evaluated the Pareto frontier of accuracy and runtime of LSQRsparse, SPIR and FOSSILS \citep{epperly2026fast}.\footnote{We tried implementing APLICUR from \citep{huh2026adaptive}, but it failed to reach the desired accuracy on the tasks considered here because the CUR approximations lead to ill-conditioned preconditioners, which caused the iterative algorithm to fail to converge within the iteration budget, so we omit it.} Across a variety of tasks with different sizes, coherence, and spectral distributions, LSQRsparse consistently achieves the same levels of accuracy $f_{\text{res}}$ (or lower) faster than the other solvers, see Figure \ref{fig::pareto_baselines}. Each curve is constructed by running the algorithms at different termination criteria, and picking the fastest one reaching any given accuracy, up to the resolution threshold defined as the point where the solutions achieved by QR and SVD disagree by more than the accuracy requested, which comes from floating point error. 

Since the solvers SPIR and FOSSILS are designed to be backward stable, they also have low forward error. Figure \ref{fig::pareto_baselines_fwd} shows that LSQRsparse also achieves small forward error much more quickly than these baselines on the tasks considered here.

\section{Hyperparameter tuning settings}\label{app::hyper}

The hyperparameter configurations considered in Section \ref{sec::analysis} are reported in Table \ref{tab::hyper}. Note that this covers a more extensive family of algorithms than the previous work on hyperparameter tuning for randomized preconditioning \citep{cho2025surrogate}. The stopping constant $c$ is set so that the iterative algorithm terminates on backward error $\textit{atol}=c\sqrt{\varepsilon}$, while the row norm/uniform mixing parameter $\alpha$ controls the tradeoff between sampling rows based on their norm or uniformly. The total search space is slightly smaller than the cartesian product of all hyperparameters, because some of them are not compatible with each other. For example, the nonzeros per row or column $k$ is only set when the chosen sketch operator is sparse sign.

The hyperparameters used during search are intended to be broad and cover a large class of algorithms, including key components in the recent randomized solvers of \citep{epperly2026fast, huh2026adaptive}, while the hyperparameters used in the refinement phase are narrowly centered around the discovered solution, i.e. $\gamma=1.9$, $\alpha=0.5$, and $c=0.16$. The refined solution selects sketch size $\gamma=2.23$, row sampling leverage $\alpha=0.58$, stopping constant $c=0.36$. To make it more competitive with LLM-guided code evolution, we also include some code-level optimization options such as which LAPACK drivers to use for certain operations.

\begin{table}[h]
    \centering
    \caption{Hyperparameter configurations used for random search and TPE on the least squares task with size $(30000,2000)$, condition number $10^{10}$, and target accuracy $10^{-4}$. }
    \label{tab::hyper}
    \small
    \setlength{\tabcolsep}{3.5pt}
    \begin{tabular}{@{}lp{0.5\linewidth}l@{}}
    \toprule
    Hyperparameter & Allowed values & Sampling \\
    \midrule
    \multicolumn{3}{@{}l}{\emph{Used during search by TPE (200 evaluations)}} \\
    \addlinespace[2pt]
    Sketch operator & sparse sign (col.), sparse sign (row), row sampling, SRHT, Gaussian & uniform \\
    Sketch size factor $\gamma$ & $[1,\ 5]$ & uniform \\
    Nonzeros per col./row $k$ & $\{1,\dots,100\}$; sparse sign only & log-uniform (int.) \\
    Row probabilities & row norms, approx.\ leverage; row sampling only & uniform \\
    Leverage mixing $\alpha$ & $[0,\ 1]$; row sampling only & uniform \\
    Replacement & without, with; row sampling only & uniform \\
    Preconditioner & QR, Cholesky, SVD, shifted CholeskyQR of $SA$ & uniform \\
    Factor precision & double, single & uniform \\
    QR routine & NumPy, SciPy, LAPACK \texttt{geqrf}; QR only & uniform \\
    Applying $R^{-1}$ & explicit inverse, SciPy or LAPACK triangular solve; not SVD & uniform \\
    Iterative solver & LSQR, LSMR, CGLS, Chebyshev, heavy ball & uniform \\
    Spectrum safety & $[0.5,\ 2]$; Chebyshev, heavy ball only & log-uniform \\
    Momentum scale & $[0,\ 1]\times$ optimal $\beta$; heavy ball only & uniform \\
    Residual replacement & never, every $16$, $32$, $64$ its.; CGLS, Chebyshev, heavy ball & uniform \\
    Warm start & none, sketch-and-solve & uniform \\
    Warm-start form & normal eq., QR of $[SA \mid Sb]$; sketch-and-solve with QR only & uniform \\
    Stopping constant $c$ & $[0.0003,\ 1]$ & log-uniform \\
    Fallback & none, LAPACK \texttt{gels} & uniform \\
    \midrule
    \multicolumn{3}{@{}l}{\emph{Used to refine the LLM-discovered programs by TPE (50 evaluations)}} \\
    \addlinespace[2pt]
    Sketch size factor $\gamma$ & $[1.2,\ 3.2]$ & uniform \\
    Leverage mixing $\alpha$ & $[0.15,\ 0.85]$ & uniform \\
    Stopping constant $c$ & $[0.05,\ 0.5]$ & log-uniform \\
    \bottomrule
    \end{tabular}
\end{table}

\section{Reproducibility of LLM-driven discovery and further comparison with hyperparameter tuning}\label{app::reproducibility}
\begin{figure}[h]
\begin{center}
    \includegraphics[width=\textwidth]{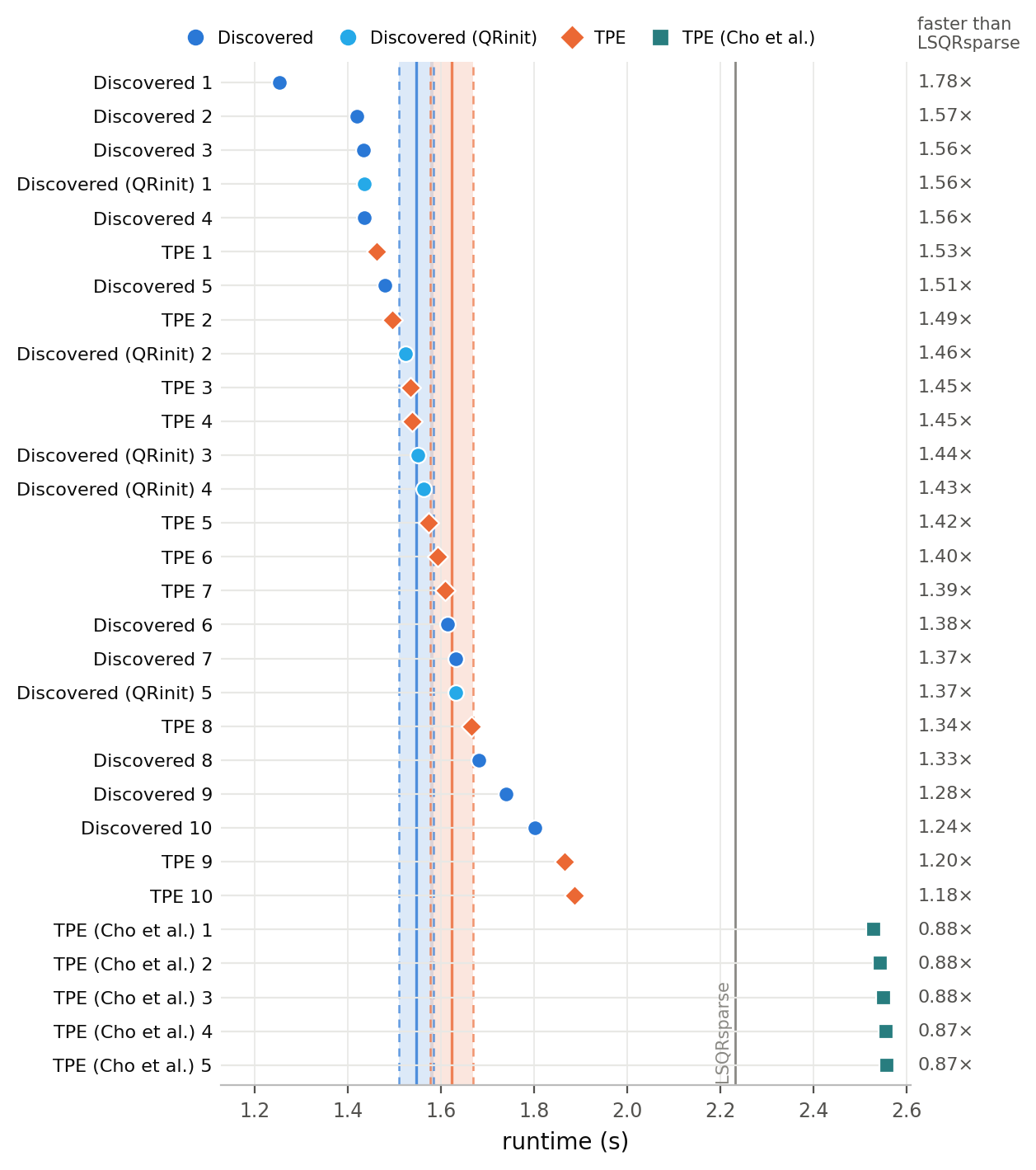}
\end{center}
\caption{Runtime achieved by the discovered runs on many independent runs of the LLM-guided discovery framework and TPE, averaged over $10$ matrices of size $30,000 \times 2000$. This includes $10$ algorithms discovered with the standard settings used throughout the paper (dark blue), $5$ algorithms discovered after starting from direct QR as initial program (light blue), $10$ algorithms discovered running TPE on the configuration described in Appendix \ref{app::hyper} (orange), and $5$ algorithms discovered running TPE on the configuration of \citep{cho2025surrogate} (teal). The mean performance of the algorithms discovered by our framework is $5\%$ faster than well-specified TPE, while miss-specified TPE is stuck at the ceiling defined by its configuration, which is much slower. The dashed lines represent the standard error of the mean, and they barely overlap (discovered: $(1.509,1.583)s$, TPE: $(1.576,1.668)s$).}
\label{fig::reproducibility+tpe}
\end{figure}

We investigated several reproducibility aspects of our LLM-guided discovery framework and hyperparameter tuning for least squares practical algorithm discovery. We focus on the task of size $(30000,2000)$, condition number $10^{10}$ and accuracy target $\varepsilon=10^{-4}$ from Section \ref{sec::analysis}.
\begin{itemize}
    \item We performed $10$ independent discovery runs using the standard settings of the paper
    \item We also investigated the effect of the initial program on the discovery, and performed $5$ independent discovery runs starting from a direct QR solver.
    \item We performed $10$ independent runs of TPE-based hyperparameter tuning with the settings of Appendix \ref{app::hyper}.
    \item We performed $5$ independent runs of TPE-based hyperparameter tuning on the grid of hyperparameters specified in \citep{cho2025surrogate}.
\end{itemize} 
Figure \ref{fig::reproducibility+tpe} summarizes the results of this experiment. The initial program does not seem to affect the final performance of the discovered algorithms, who all perform much faster than the best baseline, albeit with some variability. The well-specified TPE runs are also quite reliable in improving the baselines, but they reach inferior performance compared to LLM-guided evolution. 
On the other hand, hyperparamter tuning can only perform well if the specified grid of parameter allows it. When restricting the search to the settings studied in \citep{cho2025surrogate}, hyperparameter tuning failed to find an algorithm faster than LSQRsparse on this task.

We compared our framework with an additional baseline that uses Gaussian-Process surrogate-based autotuning, implemented through the GPSampler algorithm in Optuna.\footnote{We did not implement the exact GP surrogate in \citep{cho2025surrogate} (GPTune) because it does not allow us to use our evaluation function for the algorithm, which is essential to make the comparison robust and minimize runtime noise.}
Similar to the TPE comparison, we consider both searches over a strong configuration (Table \ref{tab::hyper}) and a weaker one (the one in \cite{cho2025surrogate}), with $5$ and $3$ independent runs, respectively. 

\begin{figure}[h]
\begin{center}
    \includegraphics[width=\textwidth]{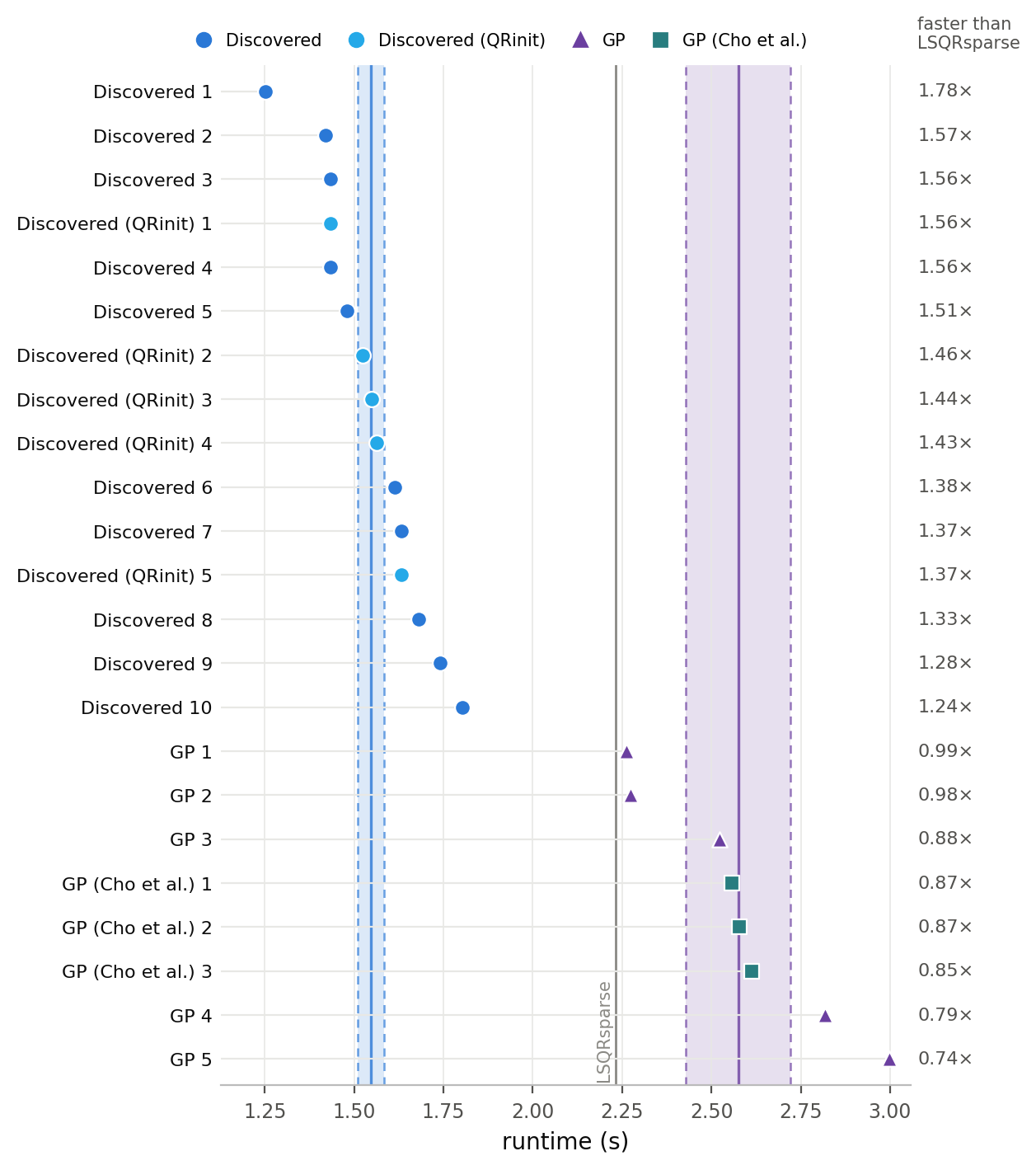}
\end{center}
\caption{Performance of algorithms discovered by our framework and by GPSampler on multiple independent runs. On average, the discovered programs are roughly $40\%$ faster than the ones found using GPSampler (the dashed lines represent the standard error of the mean).}
\label{fig::reproducibility+gp}
\end{figure}

The algorithms found on the search space of \citep{cho2025surrogate} performed similarly to the ones found by TPE on the same space, but the ones on the larger space were much slower (Figure \ref{fig::reproducibility+gp}). This is likely due to the fact that Table \ref{tab::hyper} is higher dimensional than the one of \citep{cho2025surrogate}, and contains many more categorical variables, so Bayesian optimization fails to explore the space successfully within the $200$ iteration budget. Nevertheless, the magnitude of this difference is somewhat surprising.

\section{Further analysis of discovered algorithms}\label{app::analysis}

This section contains a detailed discussion of some algorithms discovered by our framework together with the full Python programs. We will see that these algorithms can be quite different, ranging from direct to iterative methods, showing that our framework can flexibly adapt to the task at hand and discover interpretable algorithms. Moreover, the stopping criteria used in some of these algorithms are an interesting novel recombination of known principles in CGLS termination analysis.

\begin{figure}[h]
\begin{center}
    \includegraphics[width=\textwidth]{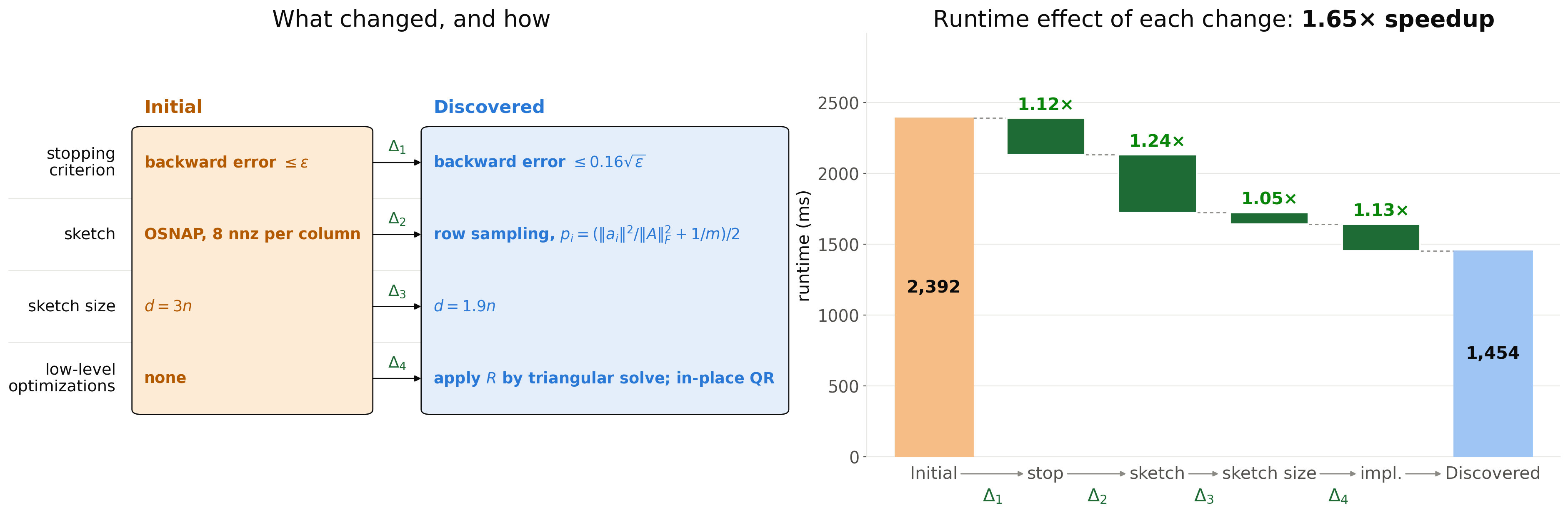}
\end{center}
\caption{Speedup decomposition of the algorithm discovered on matrices of size $(30000,2000)$, $\kappa=10^{10}$, and $\varepsilon=10^{-4}$. \textbf{Left.} Differences between the initial program (in orange) used in the discovery process and the best one (in blue). \textbf{Right.} Breakdown of performance gains from each change implemented.}
\label{fig::analysis_30k2kc10}
\end{figure}

\paragraph{More details on the task from Section \ref{sec::analysis}.} The initial program on this task achieves $f_{\text{res}} = 1+5.8 \times 10^{-8}$, well below the target of $f_{\text{res}} \le 1+10^{-4}$, so one might suspect that runtime improvements can come purely from stopping earlier. We show that the change in termination condition for the discovered algorithm only accounted for $1.12\times$ speedup, and most of the gain was due to moving away from sparse sign sketches in favor of random row sampling. See \Cref{fig::analysis_30k2kc10} for a breakdown of each improvement and (\Cref{code::disco_30k2k_1e10}, Appendix \ref{app::programs}) for the full program.

\begin{figure}[h]
\begin{center}
    \includegraphics[width=\textwidth]{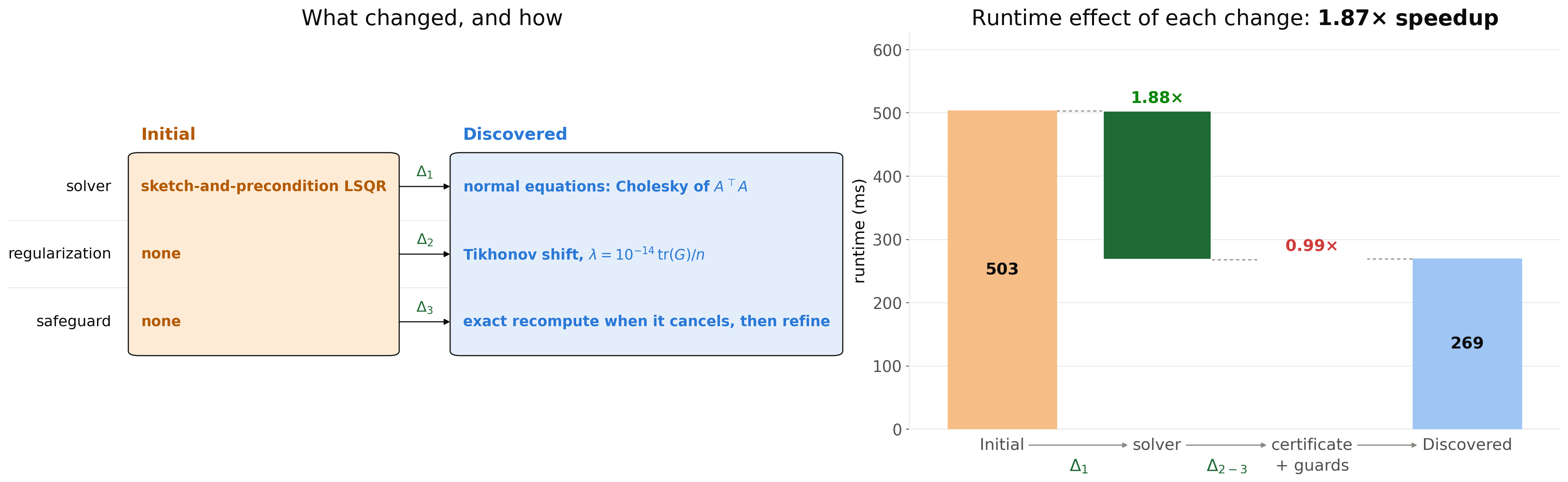}
\end{center}
\caption{Speedup decomposition of the algorithm discovered on matrices of size $(20000,1000)$, $\kappa=10^{6}$, and $\varepsilon=10^{-4}$. \textbf{Left.} Differences between the initial program (in orange) used in the discovery process and the best one (in blue). \textbf{Right.} Breakdown of performance gains from each change implemented.}
\label{fig::analysis_20k1kc6}
\end{figure}

\paragraph{Algorithm Analysis for size $(20000,1000)$ and conditioning $10^6$.} The second task is the one with size $(20000,1000)$, condition number $10^{6}$ and accuracy target $10^{-4}$, visualized in Figure \ref{fig::analysis_20k1kc6}. Here the discovered algorithm differs substantially from the initial program, switching from preconditioning to a direct Cholesky solver, which reduces runtime to roughly half. This shows that our framework can adapt to the task at hand and explore algorithms that are very different from the seed.

\paragraph{Algorithm Analysis for size $(20000,2000)$ and conditioning $10^6$.} In this task, most of the speedup came from changing the preconditioner from $R^{-1}$ coming from the QR decomposition of the sketched data matrix to $L^{-\top}$ coming from the Cholesky decomposition of the sketched Gram matrix. 

An interesting discovery in this task came from the termination condition. The accuracy target on \Cref{eq::acc_metric} is equivalent to 
\begin{equation}\label{eq::target_fwd}
    \|A(x-x^*)\|_2^2 \le \left(1 - \frac{1}{(1+\varepsilon)^2} \right) \|Ax - b\|_2^2 = \eta \|Ax - b\|_2^2
\end{equation}
and the left-hand side can be written in terms of the Gram matrix and the gradient of the quadratic objective function $g=A^\top (Ax - b)$ as 
\begin{equation*}
    \|A(x-x^*)\|_2^2 = g^\top (A^\top A)^{-1} g
\end{equation*}
Since the sketched matrix $A_s$ contains a subset of the rows of $A$, we have $A_s^\top A_s \preceq A^\top A$, so we can upper bound $f_\text{res}$ exactly through
\begin{equation*}
    g^\top (A_s^\top A_s)^{-1} g \le 0.9\eta \|Ax - b\|_2^2
\end{equation*}
This uses only information available to the iterative algorithm, i.e. the sketched Gram matrix and the current residual. The $0.9$ factor is added for slack. Although mathematically well-motivated, this stopping criterion does not provide a speedup on this task.

\begin{figure}[h]
\begin{center}
    \includegraphics[width=\textwidth]{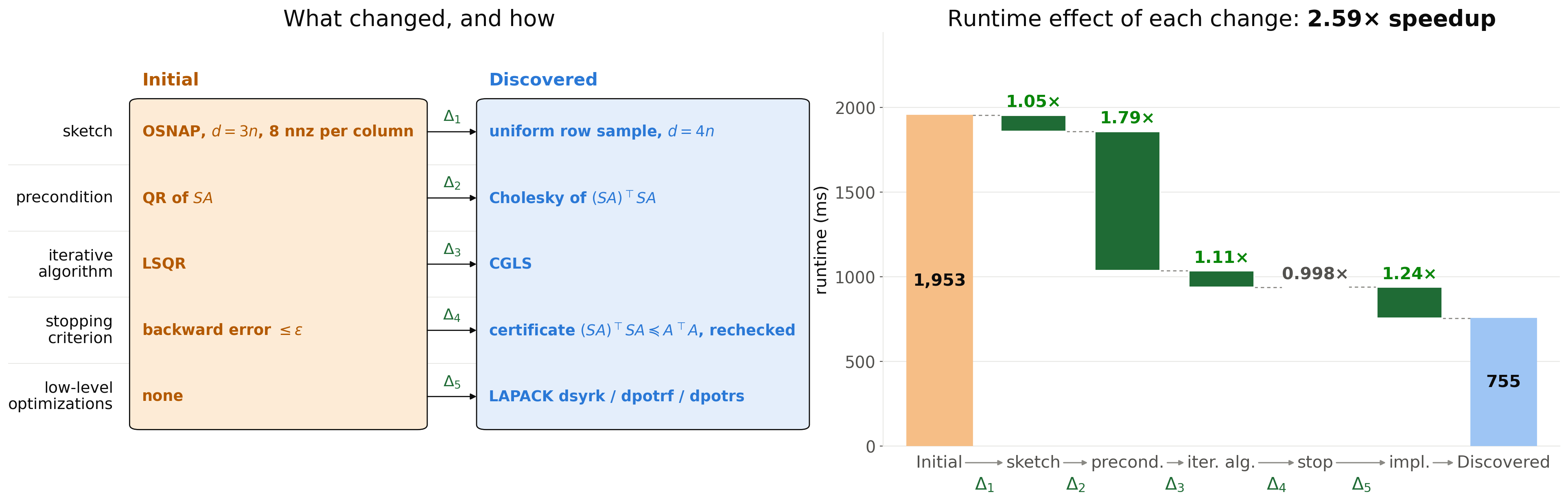}
\end{center}
\caption{Speedup decomposition of the algorithm discovered on matrices of size $(20000,2000)$, $\kappa=10^{6}$, and $\varepsilon=10^{-4}$. \textbf{Left.} Differences between the initial program (in orange) used in the discovery process and the best one (in blue). \textbf{Right.} Breakdown of performance gains from each change implemented.}
\label{fig::analysis_20k2kc6}
\end{figure}

\begin{figure}[h]
\begin{center}
    \includegraphics[width=\textwidth]{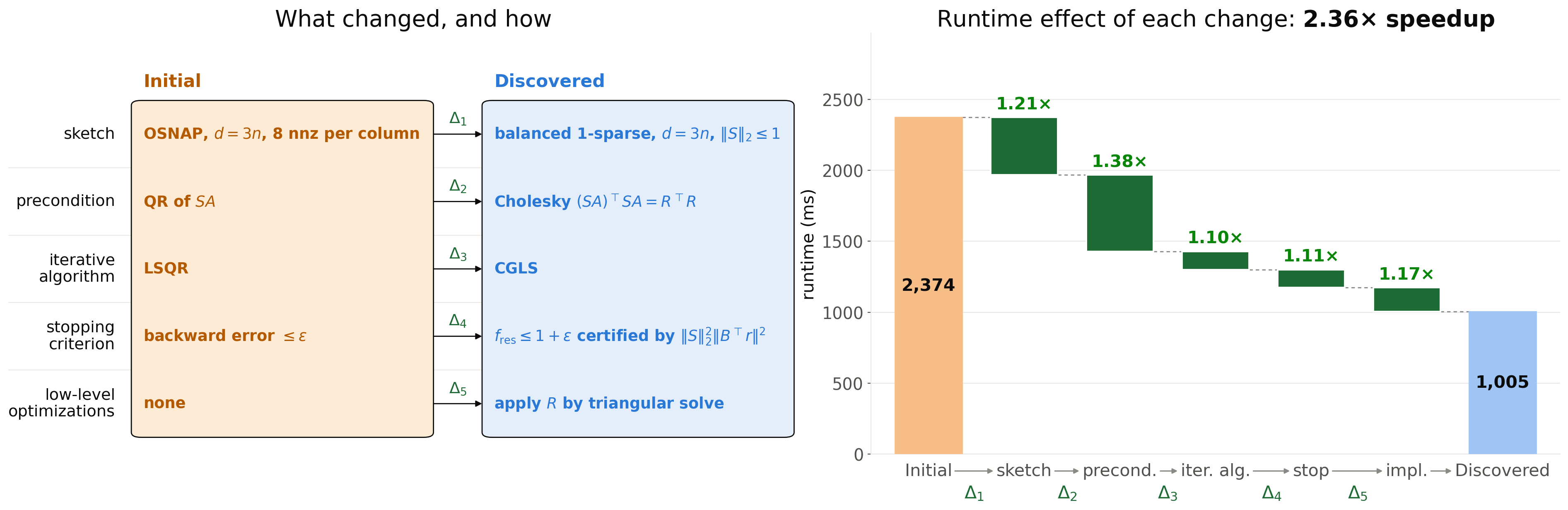}
\end{center}
\caption{Speedup decomposition of the algorithm discovered on matrices of size $(30000,2000)$, $\kappa=10^{6}$, and $\varepsilon=10^{-4}$. \textbf{Left.} Differences between the initial program (in orange) used in the discovery process and the best one (in blue). \textbf{Right.} Breakdown of performance gains from each change implemented.}
\label{fig::analysis_30k2kc6}
\end{figure}

\paragraph{Algorithm Analysis for size $(30000,2000)$ and conditioning $10^6$.} The algorithm discovered for this task was also a randomized iterative solver, but with a $1$-sparse sketch followed by a Cholesky-preconditioned CGLS. The termination condition is also non-standard, and we can explain it as follows. Let $r = Ax - b$, and $B = A R^{-1}$, where $R$ is the Cholesky factor of $(SA)^\top (SA)$, i.e. the preconditioner. The discovered stopping criterion is 
\begin{equation}\label{eq::newbound2}
    \|B^\top r\|_2 \le \frac{\sqrt{\eta}}{1.1 \|S\|_2} \|r\|_2
\end{equation}
where $\eta = \left( 1 - \frac{1}{(1+\varepsilon)^2} \right)$. This is a provably sufficient stopping criterion for the accuracy target in exact arithmetic. Starting from the equivalent condition \Cref{eq::target_fwd} for $f_{\text{res}} \le 1+\varepsilon$, the right-hand side can be written as 
\begin{equation}\label{eq::spectral_bound}
    \|A(x - x^*)\|_2^2 = (B^\top r)^\top (B^\top B)^{-1} (B^\top r) \le \frac{\|B^\top r\|_2^2}{\sigma_{\min} (B)^2}
\end{equation}
Note that the matrix $SB$ has orthonormal columns, because $B=AR^{-1}$ where $R$ is the Cholesky factor of $(SA)^\top (SA)$. Therefore, for any vector $y \in \mathbb{R}^n$,
\begin{align*}
    \|y\|_2 = \|SBy\|_2 \le \|S\|_2 \|By\|_2
\end{align*}
This implies that $\sigma_{\min}(B) \ge \|S\|_2^{-1}$. Combining this with \Cref{eq::spectral_bound} and \Cref{eq::target_fwd} shows that the stopping criterion implies $f_{\text{res}} \le 1+\varepsilon$. The algorithm builds a sketch matrix with (normalized) $1$-sparse columns and rows with disjoint support, so they are orthonormal, hence $\|S\|_2=1$ by construction and does not need to be computed. This is what makes the bound \Cref{eq::newbound2} efficiently computable by the iterative algorithm, since $\|r\|$ and $\|B^\top r\|$ are already available by CGLS iterates. Many of these ingredients are known in standard CGLS literature \citep{bjorck2024numerical, saad2003iterative}, and related spectral analysis of preconditioners appears in \citep{arioli2015preconditioning}. However, these ideas have not been made practical in standard scientific computing libraries, and, to the best of our knowledge, the application to $f_{\text{res}}$ instead of backward error is novel. Contrary to the previous task, the new stopping criterion provided a $1.11\times$ runtime speedup relative to the default CGLS stopping rule (\Cref{fig::analysis_30k2kc6}).

\begin{figure}[h]
\begin{center}
    \includegraphics[width=\textwidth]{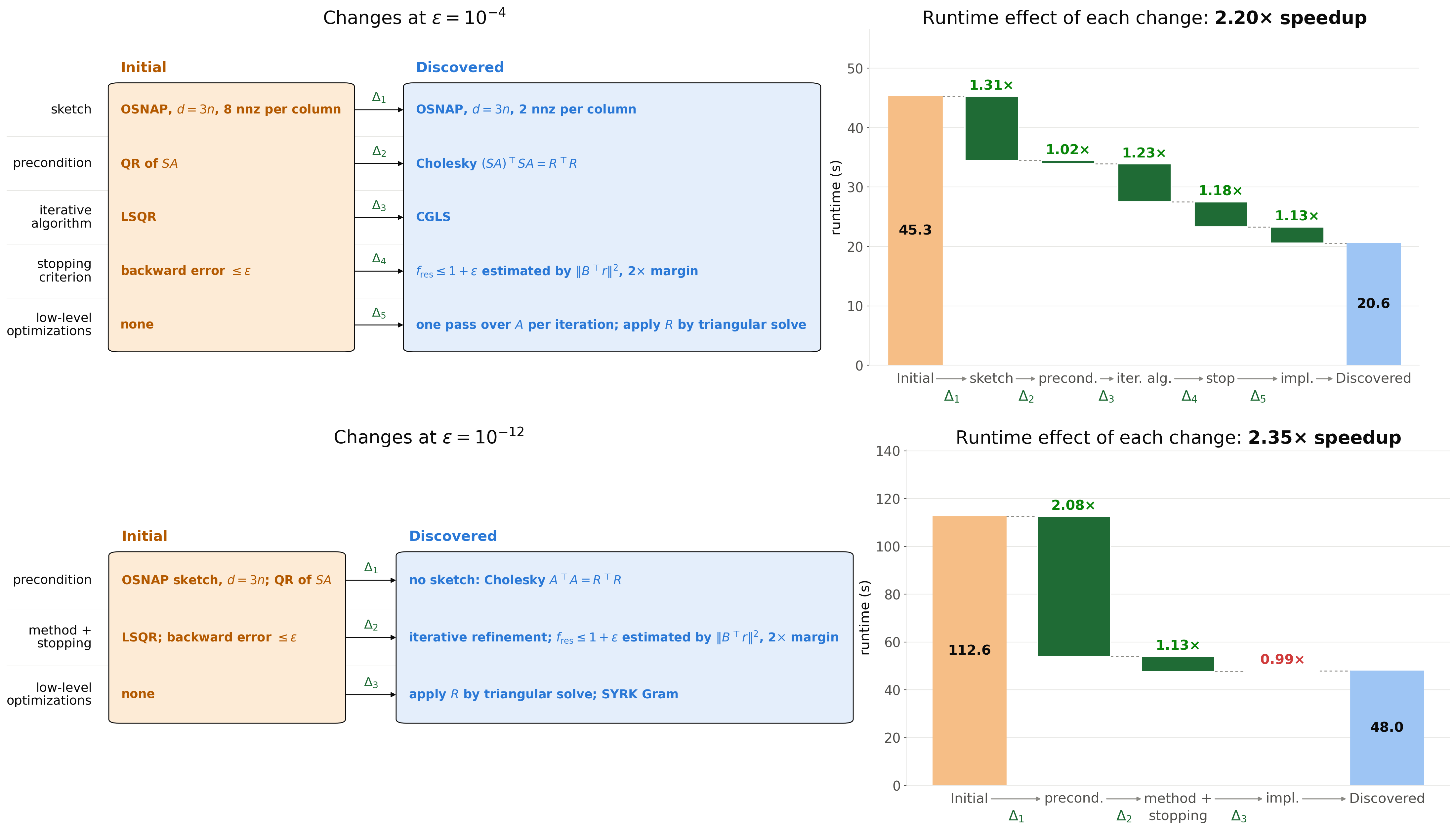}
\end{center}
\caption{\textbf{Top.} Changes introduced by the discovery process for the SUSY task at $\varepsilon=10^{-4}$, and their contribution to the speedup. \textbf{Bottom.} The same decomposition for $\varepsilon=10^{-12}$.}
\label{fig::analysis_susy}
\end{figure}

\paragraph{Algorithm Analysis for the SUSY task}
Figure \ref{fig::analysis_susy} presents the runtime gain decomposition for the algorithms discovered on the SUSY task at medium accuracy $\varepsilon=10^{-4}$ and high accuracy $\varepsilon=10^{-12}$. The former is a randomized preconditioner based OSNAP and CGLS, while the latter is a direct Cholesky solve with iterative refinement.

\section{More details on the SUSY task}\label{app::susy_deepdive}
\begin{figure}[h]
\begin{center}
    \includegraphics[width=\textwidth]{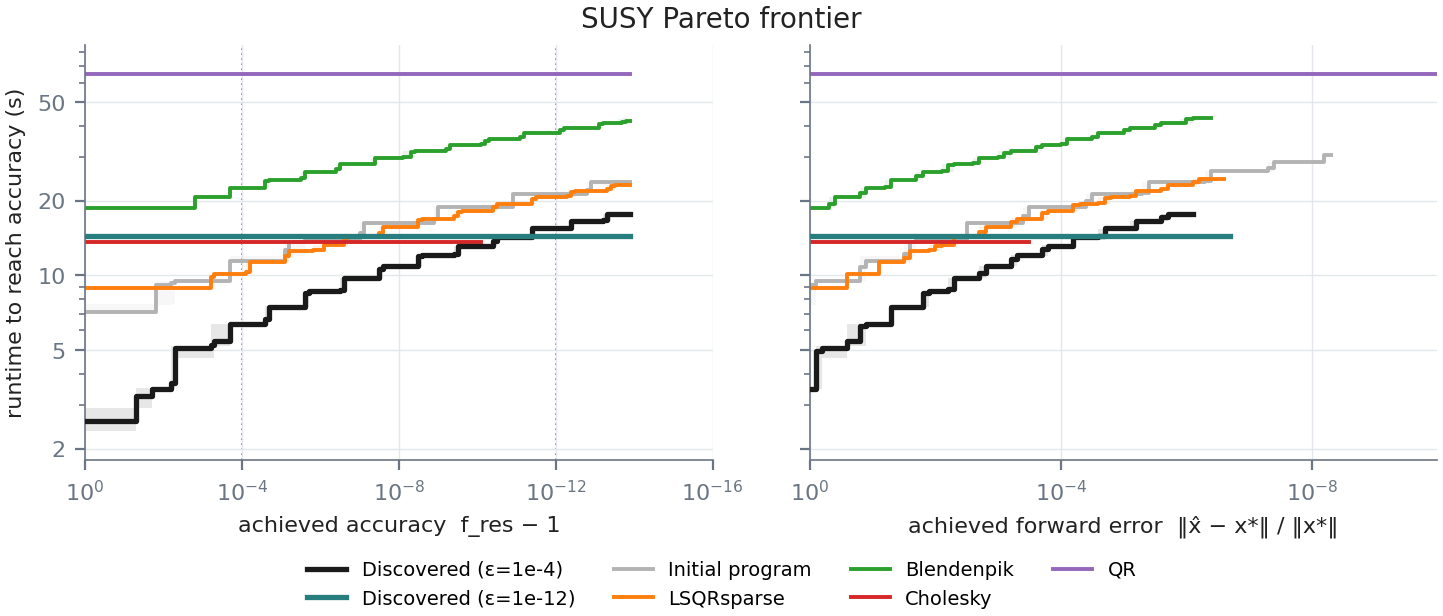}
\end{center}
\caption{Runtime vs relative residual (and forward error) Pareto frontier for the SUSY task. The algorithm discovered for a medium accuracy target is Pareto optimal for any $\varepsilon \ge 10^{-11}$, while the one discovered for a high-accuracy target dominates the end of the curve.}
\label{fig::pareto_susy}
\end{figure}

In this section, we provide more information on the performance of the discovered algorithms for the SUSY task, in particular exploring more fine-grained instance-level statistics and the Pareto frontier of accuracy and runtime.
\paragraph{Instance-level statistics} We focus on the SUSY task at target accuracy $\varepsilon = 10^{-4}$. To obtain more test instances, we split the test dataset in ten matrices with $300,000$ rows instead of three matrices with $1,000,000$ rows. Table \ref{tab::susy_instances} shows that the performance of the discovered algorithm is very consistent across instances. This is caused by concentration of measure, since at these scales many important quantities influencing the performance of randomized iterative solvers like the leading singular values concentrate strongly, even though the entries of the matrices can be quite different. Each instance is run $10$ times on different random seeds, and we report the key quantiles of the runtime distribution, together with maximum relative residual error achieved.
\begin{table}[t]
\centering
\small
\setlength{\tabcolsep}{4pt}
\caption{Per-instance robustness of the discovered program on SUSY at $\varepsilon = 10^{-4}$.
Ten disjoint instances of $300{,}000\times2000$ (Gaussian-kernel least squares on SUSY), cut from rows disjoint from those used during discovery. }
\label{tab::susy_instances}
\begin{tabular}{@{}c rrr rr r r@{}}
\toprule
 & \multicolumn{3}{c}{time (s)} & & & & \\
\cmidrule(lr){2-4}
Instance & median & q10--q90 & min--max & max $f_{\text{res}}-1$ & violations & fwd.\ error & speedup vs LSQRsparse \\
\midrule
0 & 6.35 & 6.33--6.36 & 6.33--6.38 & $2.7\times10^{-5}$ & 0/10 & $5.8\times10^{-2}$ & $1.98\times$ \\
1 & 6.35 & 5.78--6.36 & 5.78--6.36 & $6.7\times10^{-5}$ & 0/10 & $5.7\times10^{-2}$ & $1.99\times$ \\
2 & 6.34 & 6.31--6.35 & 6.31--6.38 & $2.8\times10^{-5}$ & 0/10 & $5.6\times10^{-2}$ & $1.99\times$ \\
3 & 6.35 & 6.34--6.40 & 6.34--6.55 & $2.6\times10^{-5}$ & 0/10 & $5.1\times10^{-2}$ & $1.99\times$ \\
4 & 6.33 & 6.33--6.35 & 6.32--6.36 & $2.8\times10^{-5}$ & 0/10 & $5.2\times10^{-2}$ & $1.99\times$ \\
5 & 6.35 & 6.29--6.36 & 5.78--6.36 & $5.9\times10^{-5}$ & 0/10 & $5.8\times10^{-2}$ & $1.99\times$ \\
6 & 6.34 & 6.32--6.35 & 6.30--6.36 & $3.0\times10^{-5}$ & 0/10 & $5.8\times10^{-2}$ & $1.99\times$ \\
7 & 6.36 & 6.35--6.39 & 6.34--6.40 & $2.7\times10^{-5}$ & 0/10 & $5.3\times10^{-2}$ & $1.99\times$ \\
8 & 6.34 & 5.77--6.36 & 5.77--6.39 & $6.6\times10^{-5}$ & 0/10 & $6.1\times10^{-2}$ & $1.98\times$ \\
9 & 6.35 & 6.34--6.36 & 6.34--6.37 & $2.8\times10^{-5}$ & 0/10 & $5.1\times10^{-2}$ & $1.94\times$ \\
\midrule
All & $6.35 \pm 0.00$ & & & $6.7\times10^{-5}$ & 0/100 & $5.5\times10^{-2}$ & $1.98\times$ {\scriptsize [1.97, 1.99]} \\
\bottomrule
\end{tabular}
\end{table}
\paragraph{Speed vs accuracy Pareto frontier} Even though the algorithms discussed in Table \ref{tab::susy} achieve a lower accuracy than the baselines on the medium accuracy task $\varepsilon=10^{-4}$, the runtime improvement still holds after accounting for the exact accuracy achieved by the algorithms. Figure \ref{fig::pareto_susy} shows the Pareto frontier for both relative residual and forward error. In terms of relative residual, the discovered algorithm improvement over the baselines can be as large as $3.5\times$ at low accuracy ($\varepsilon=0.1$). This figure also uses the $10$ smaller matrices of size $300,000 \times 2000$, so we can generate error bars using the maximum and minimum runtime, but those are again barely visible.


\section{Forward and backward error of discovered algorithms}\label{app::fwdbwd}
\begin{figure}[h]
\begin{center}
    \includegraphics[width=\textwidth]{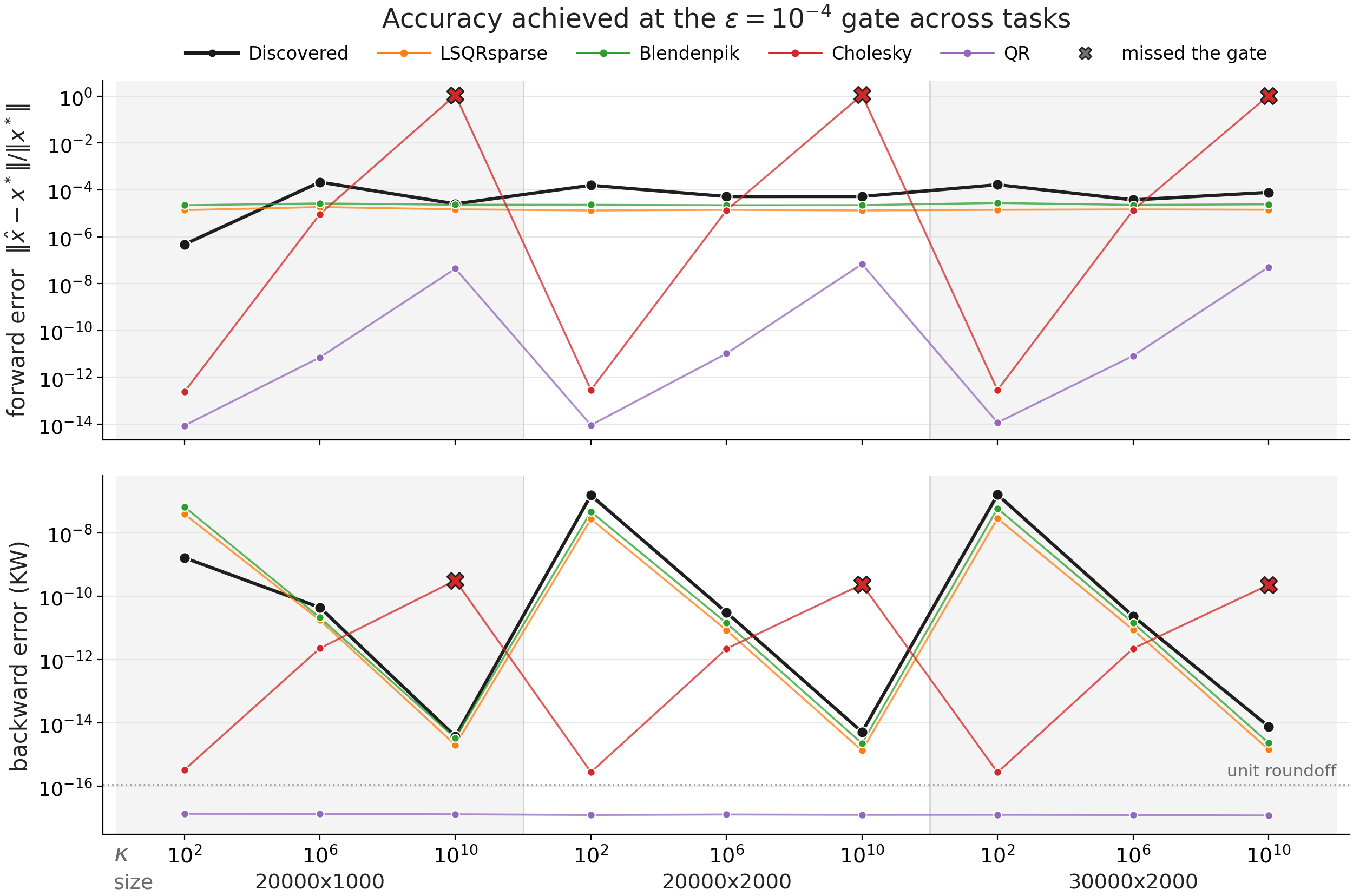}
\end{center}
\caption{\textbf{Top.} Forward error achieved by different algorithms on the SUSY task. Each datapoint is the average runtime over the fastest configuration of the algorithm reaching the target error over a test set of $10$ random matrices. \textbf{Bottom.} Same as above for the KW estimate of the backward error.}
\label{fig::fwdbwd}
\end{figure}

The algorithms discovered in the main paper are designed to be fast conditioned on achieving the relative residual accuracy specified by the task. This does not directly control other important accuracy metrics such as forward and backward errors. Recall that the forward error of a vector $x$ is defined as $\|x^* - x\|_2 / \|x^*\|_2$, while the backward error is the norm of the perturbation to the input data required to produce an instance for which $x$ is an exact solution. This is difficult to compute exactly, and a standard alternative that is easier to compute is the Karlson-W\'alden estimator \citep{karlson1997estimation}. Figure \ref{fig::fwdbwd} shows the achieved forward and backward errors of the algorithms discussed in Section \ref{sec::experiments}, Figure \ref{fig::4.2}. The discovered algorithms have comparable forward and backward errors to the baselines belonging to the same class of classical randomized iterative solvers. Table \ref{tab:susy_fwdbwd} shows the forward and backward error of the discovered algorithms and baselines for the SUSY task, which follow a similar pattern. 
If the focus of the analysis changes to prioritizing these errors, the baselines should also change to include \cite{epperly2026fast}, but it is beyond the scope of this paper.
\begin{table}[b]
\centering
\caption{Forward and backward error on the SUSY task.}
\label{tab:susy_fwdbwd}
\begin{tabular}{l rr rr}
\toprule
 & \multicolumn{2}{c}{$\varepsilon = 10^{-4}$} & \multicolumn{2}{c}{$\varepsilon = 10^{-12}$} \\
\cmidrule(lr){2-3} \cmidrule(lr){4-5}
Program & forward & backward & forward & backward \\
\midrule
Discovered & $6.2\times10^{-2}$ & $2.9\times10^{-8}$ & $2.1\times10^{-7}$ & $2.3\times10^{-15}$ \\
LSQRsparse & $1.3\times10^{-2}$ & $5.9\times10^{-9}$ & $2.8\times10^{-6}$ & $1.2\times10^{-12}$ \\
Blendenpik & $2.3\times10^{-2}$ & $1.1\times10^{-8}$ & $4.1\times10^{-6}$ & $1.9\times10^{-12}$ \\
Cholesky & $3.8\times10^{-4}$ & $5.0\times10^{-12}$ & \textcolor[HTML]{C44E52}{$3.8\times10^{-4}$} & \textcolor[HTML]{C44E52}{$5.0\times10^{-12}$} \\
QR & $3.3\times10^{-11}$ & $7.7\times10^{-17}$ & $3.3\times10^{-11}$ & $7.7\times10^{-17}$ \\
\bottomrule
\end{tabular}
\end{table}
\section{(Even more) out-of-distribution generalization}\label{app::OOD}
In this section we investigate the behavior of the discovered algorithms on tasks that differ from the training settings more strongly than in Section \ref{sec::analysis}. We focus on the algorithm discovered for the synthetic task with size $(30000,2000)$, target accuracy $\varepsilon=10^{-4}$, smoothly decaying singular values, condition number $10^{10}$, uniform leverage scores and optimal residual $r^*=0.01$. This is the algorithm discussed in Figure \ref{fig::breakdown_vs_hyper}, and we chose it because it performs row sampling, which is known to be brittle when the leverage scores of the rows concentrate on a few rows, which coincides with having high coherence.
We varied the following characteristics of the task.
\begin{itemize}
    \item \textbf{Spectral distribution}: Instead of smoothly decaying singular values, consider clustered singular value distributions, where the singular values are tightly clustered around $1$ and $10^{-10}$ in equal amounts (\textit{clust2}) or around $1, 10^{-3}, 10^{-5},10^{-8}$ and $10^{-10}$ in equal amounts (\textit{clust5}).
    \item \textbf{Optimal residual}: Consider both optimal residual $r^*=0.01$ and less consistent instances with $r^*=0.5$.
    \item \textbf{Coherence level}: We consider both tasks with unit coherence (uniform leverage) or high coherence (skewed leverage), following the T1 construction from \citep{ma14statistical}.
\end{itemize}
Combining all these features leads to $8$ tasks, and the performance of the discovered algorithm on them is collected in Table \ref{tab::ood_disco}. The improvement over the best baseline generalizes across all tasks except for the high-coherence ones, where the algorithm fails to solve the least squares instances to any useful level of accuracy. This is not surprising, because these instances are very different from the task provided to our discovery framework. Different instances from the typical ones are likely to arise in practice, but hopefully the task of interest is well-defined enough to observe a few outlier instances in the training set. To simulate this, we discover an algorithm for a task which contains $8$ matrices from the original training set, and $2$ matrices with high coherence. When tested on only high-coherence instances, this algorithm still fails to reach the target accuracy on a few instances, but it is much more robust than the original one.

\begin{table*}[!t]
\centering
\small
\setlength{\tabcolsep}{5pt}
\caption{Out-of-distribution test of two discovered programs at the Fig.~4 setting
($30000\times2000$, $\kappa=10^{10}$, $\varepsilon=10^{-4}$). \textbf{Disco1} is the original program (Figure \ref{fig::breakdown_vs_hyper}). \textbf{Disco2} was discovered by the same search with 2 of those 10 replaced by high-coherence matrices. Each number is the result of timing 5 matrices on 3 different seeds, for a total of 15 program evaluations. The instances on which \textbf{Disco2} fails the accuracy target of $\varepsilon=10^{-4}$ have error at most $2.7 \times 10^{-4}$.}
\label{tab::ood_disco}
\begin{tabular}{@{}l r rrr rrr@{}}
\toprule
 & LSQRsparse & \multicolumn{3}{c}{Disco1} & \multicolumn{3}{c}{Disco2} \\
\cmidrule(lr){3-5}\cmidrule(lr){6-8}
Task & Time & Pass & Time & Max Err. & Pass & Time & Max Err. \\
\midrule
smooth, unif, resid0.01 & 2255 & 15 & 1440 & 4.5e-5 & 15 & 1522 & 6.0e-5 \\
\midrule
clust2, unif, resid0.01 & 2255 & 15 & 1449 & 4.2e-5 & 15 & 1528 & 6.7e-5 \\
clust2, unif, resid0.5  & 2278 & 15 & 1465 & 4.5e-5 & 15 & 1544 & 6.8e-5 \\
clust5, unif, resid0.01 & 2248 & 15 & 1443 & 4.9e-5 & 15 & 1525 & 6.5e-5 \\
clust5, unif, resid0.5  & 2237 & 15 & 1436 & 4.5e-5 & 15 & 1519 & 6.3e-5 \\
\midrule
clust2, high, resid0.01 & 2273 & \textcolor{red}{0} & 3413 & 4.0  & \textcolor{orange}{13} & 1679 & 1.4e-4 \\
clust2, high, resid0.5  & 2320 & \textcolor{red}{0} & 2863 & 2.4  & \textcolor{orange}{13} & 1719 & 1.1e-4 \\
clust5, high, resid0.01 & 2223 & \textcolor{red}{0} & 3364 & 0.83 & \textcolor{orange}{11} & 1665 & 2.7e-4 \\
clust5, high, resid0.5  & 2237 & \textcolor{red}{0} & 2644 & 23   & \textcolor{orange}{13} & 1649 & 1.1e-4 \\
\bottomrule
\end{tabular}
\end{table*}

\section{Other ideas explored}
This appendix collects a couple of approaches that were considered for part of this project and later dropped from the final version of our framework, but they might be of interest to some readers.

\subsection{Deterministic robust proxy for runtime}
Since measured runtime is noisy and depends on many hardware and software-specific factors that are hard to control, we investigated using a deterministic proxy for runtime as score function for ShinkaEvolve. We took inspiration from the Roofline Model \citep{williams2009roofline, gholami2024ai} to measure the ``speed'' of each operation according to its arithmetic intensity. Given a machine with peak compute throughput $P$ (i.e. number of FLOPs per second) and peak memory throughput $W$ (i.e. number of bytes moved per second), the performance of an operation $op$ under the roofline model is given by
\begin{equation*}
    P_{op} = \min \begin{cases}
        P \\
        W \cdot I_{op}
    \end{cases}
\end{equation*}
where $I_{op}$ is the arithmetic intensity, i.e. the ratio between the number of FLOPs performed and the number of bytes moved. Since the peak compute throughput is typically much larger than the peak memory throughput, having higher arithmetic intensity leads to higher performance. We calibrated the throughput of general matrix-matrix multiplication on our machine and defined the weight of an operation as 
\begin{equation*}
    w(op) = \frac{P_{gemm}}{P_{op}}
\end{equation*}
We then used it to build a ``wrapper'' library around Numpy and Scipy containing common linear algebra operations used by least squares solvers, which also returned the FLOPs performed weighted by their throughput. For a program $p$ composed of basic linear algebra operations $op$, its weighted cost is
\begin{equation*}
    C_{p} = \sum_{op \in p} \text{FLOPS}(op)w(op)
\end{equation*}
Rather than reproducing the precise runtime, the purpose of the proxy is to rank programs according to their runtime, so it can be used as score function during the discovery phase. This was especially promising for least squares algorithms, which use a composition of BLAS-2 and BLAS-3 operations. This succeeded in reducing the variance of runtime measurement, however it introduced too much bias which made it much less competitive against the final version of the score function we use. There are two main issues: (i) a single roofline is often not sufficient to model performance on modern machines with different memory levels and accelerators, and (ii) having to describe the library to the LLM prompt introduced a lot of bias in the proposed algorithms.

\subsection{Minimizing residual error subject to a time budget}
Our evaluation function is designed to minimize runtime subject to a target accuracy level. The reversed question is also interesting, i.e. finding the algorithm that maximizes accuracy subject to a time budget. We ran a few experiments with this score function and the algorithms we discovered were often inferior to the ones found by our current setup. One reason is that when the same accuracy can be reached by two algorithms that are both inside the time budget, this objective function does not discriminate between them and might choose the slower one. 

\section{Discovered Practical Programs for Least Squares}\label{app::programs}
In this section, we present the complete programs discovered by our framework for the tasks discussed in Appendix \ref{app::analysis}. The full list of programs for all tasks considered in the paper will be released in a public repository at the time of publication.

\begin{lstlisting}[
  language=Python,
  caption={Program discovered for the task of size $(30000,2000)$, conditioning $10^{10}$, and accuracy target $10^{-4}$},
  label={code::disco_30k2k_1e10},
  basicstyle=\ttfamily\tiny,  % Change font size here
]
# EVOLVE-BLOCK-START
import numpy as np
from scipy import linalg as sla
from scipy.linalg.lapack import dtrtrs
from scipy.sparse.linalg import LinearOperator, lsqr

_TINY2 = 1e-40

def solve(A, b, eps, iter_budget):
    """Solve min_x ||A x - b||_2 using sampled-QR preconditioned LSQR."""
    m, n = A.shape

    if n == 0:
        return np.empty(0, dtype=np.float64)

    # A modest reduction from 2n lowers both sampled-matrix construction and
    # QR cost while retaining a well-conditioned preconditioned operator.
    d = min(m, max(n + 16, int(1.9 * n)))

    if d == m:
        SA = np.array(A, dtype=np.float64, order="F", copy=True)
        Sb = np.array(b, dtype=np.float64, copy=True)
    else:
        # Reuse the energy array as the probability array to reduce allocations
        # and full memory passes.
        probabilities = np.einsum("ij,ij->i", A, A, optimize=True)
        total = float(probabilities.sum())

        if not np.isfinite(total) or total <= 0.0:
            return np.zeros(n, dtype=np.float64)

        probabilities *= 0.5 / total
        probabilities += 0.5 / m
        probabilities /= probabilities.sum()

        # Sampling without replacement avoids wasting sketch rows on duplicates.
        idx = np.random.choice(m, size=d, replace=False, p=probabilities)

        inclusion = d * probabilities[idx]
        np.minimum(inclusion, 1.0, out=inclusion)
        scale = 1.0 / np.sqrt(inclusion)

        SA = np.array(A[idx], dtype=np.float64, order="F", copy=True)
        SA *= scale[:, None]
        Sb = b[idx] * scale

    # Form this before QR overwrites the sampled matrix.
    rhs = SA.T @ Sb

    R = sla.qr(
        SA,
        mode="r",
        overwrite_a=True,
        check_finite=False
    )
    R = np.asfortranarray(R)

    # Warm start from the sampled normal equations.
    w, info = dtrtrs(R, rhs, lower=0, trans=1, overwrite_b=1)
    if info != 0:
        w = sla.solve_triangular(
            R, rhs, trans="T", lower=False, check_finite=False
        )
        
    x0, info = dtrtrs(R, w, lower=0, trans=0, overwrite_b=1)
    if info != 0:
        x0 = sla.solve_triangular(
            R, w, lower=False, check_finite=False
        )

    residual = b - A @ x0
    residual_norm2 = float(residual @ residual)
    if residual_norm2 <= _TINY2:
        return x0

    AT = A.T

    def matvec(v):
        z, _ = dtrtrs(R, v, lower=0, trans=0)
        return A @ z

    def rmatvec(u):
        g = AT @ u
        z, _ = dtrtrs(R, g, lower=0, trans=1, overwrite_b=1)
        return z

    operator = LinearOperator(
        (m, n),
        matvec=matvec,
        rmatvec=rmatvec,
        dtype=np.float64
    )

    # Objective error is quadratic in the normalized normal-equation residual
    # after sketch preconditioning. The cap also prevents an excessively loose
    # test when callers request low accuracy.
    target = max(float(eps), 0.0)
    atol = min(1.0e-2, 0.16 * np.sqrt(target))

    y = lsqr(
        operator,
        residual,
        atol=atol,
        btol=0.0,
        conlim=0.0,
        iter_lim=int(iter_budget),
        show=False
    )[0]

    correction, info = dtrtrs(
        R, y, lower=0, trans=0, overwrite_b=1
    )
    if info != 0:
        correction = sla.solve_triangular(
            R, y, lower=False, check_finite=False
        )

    return x0 + correction
# EVOLVE-BLOCK-END
# ...
\end{lstlisting}

\begin{lstlisting}[
  language=Python,
  caption={Program discovered for the task of size $(20000,1000)$, conditioning $10^{6}$, and accuracy target $10^{-4}$},
  label={code::disco_20k1k_1e6},
  basicstyle=\ttfamily\tiny,  
]
# EVOLVE-BLOCK-START
import numpy as np
from scipy import linalg as sla

def solve(A, b, eps, iter_budget):
    """
    Super-fast overdetermined least-squares solver:
    1) Solves normal equations via fast Cholesky factorization.
    2) Estimates optimal residual and gradient in O(n^2) to perform
       the exact stopping check without doing costly O(mn) matvecs.
    3) Falls back to exact computation and one-step refinement if needed.
    """
    m, n = A.shape

    # 1. Compute AtA and Atb using optimized NumPy GEMM/GEMV
    AtA = A.T @ A
    Atb = A.T @ b

    # 2. Diagonal boost for stability
    tr = float(AtA.diagonal().sum())
    lam = max(1e-14 * tr / n, 1e-20)

    # Add the boost to AtA in-place (AtA is a temporary copy, so this is safe)
    np.fill_diagonal(AtA, AtA.diagonal() + lam)

    # 3. Cholesky factorization
    try:
        c, lower = sla.cho_factor(AtA, lower=False, overwrite_a=False, check_finite=False)
        x = sla.cho_solve((c, lower), Atb, overwrite_b=False, check_finite=False)
    except Exception:
        # Robust fallback for extreme cases
        return sla.lstsq(A, b, lapack_driver="gelsy", check_finite=False)[0]

    # 4. Ultra-fast exact stopping test
    b_norm2 = b.dot(b)
    if b_norm2 <= 1e-20:
        return np.zeros(n, dtype=A.dtype)

    tol_factor = (2.0 * eps + eps**2) / (1.0 + eps)**2

    # O(n^2) gradient and residual estimation
    AtA_x = AtA @ x
    unreg_AtA_x = AtA_x - lam * x
    g = Atb - unreg_AtA_x
    r_norm2 = b_norm2 - 2.0 * x.dot(Atb) + x.dot(unreg_AtA_x)

    # If the residual might have had catastrophic cancellation (highly consistent system),
    # fall back to exact O(mn) calculation.
    if r_norm2 < 1e-11 * b_norm2:
        r = b - A @ x
        r_norm2 = r.dot(r)
        g = A.T @ r

    dx = sla.cho_solve((c, lower), g, check_finite=False)
    w_norm2 = max(g.dot(dx), 0.0)

    # Typical case: returns immediately here after almost zero extra cost
    if w_norm2 <= tol_factor * r_norm2:
        return x

    # 5. Iterative Refinement
    x += dx

    # Small budget loop for refinement; usually 1 step is more than enough
    for _ in range(min(5, iter_budget)):
        r = b - A @ x
        r_norm2 = r.dot(r)
        if r_norm2 <= 1e-20:
            return x
        g = A.T @ r
        dx = sla.cho_solve((c, lower), g, check_finite=False)
        w_norm2 = max(g.dot(dx), 0.0)
        if w_norm2 <= tol_factor * r_norm2:
            return x
        x += dx

    return x
# EVOLVE-BLOCK-END
# ...
\end{lstlisting}

\begin{lstlisting}[
  language=Python,
  caption={Program discovered for the task of size $(20000,2000)$, conditioning $10^{6}$, and accuracy target $10^{-4}$},
  label={code::disco_20k2k_1e6},
  basicstyle=\ttfamily\tiny,  
]
# EVOLVE-BLOCK-START
import numpy as np
from scipy import linalg as sla
from scipy.linalg.blas import dsyrk
from scipy.linalg.lapack import dpotrf, dpotrs, dposv

def solve(A, b, eps, iter_budget):
    """
    Row-subsampled Cholesky-preconditioned CG for dense least squares.

    For g = A.T @ (b - A@x), the excess squared residual is
        g.T @ inv(A.T@A) @ g.
    If G = As.T@As is formed from any subset of A's rows, then
        G <= A.T@A
    in the positive-semidefinite order, and consequently
        g.T @ inv(A.T@A) @ g <= g.T @ inv(G) @ g.
    This supplies a cheap, observable stopping certificate.
    """
    m, n = A.shape
    At = A.T
    requested = max(0, int(iter_budget))
    tolerance = float(eps)

    # Four rows per unknown generally gives a strong preconditioner while
    # requiring only 40% of the full Gram work when m/n is ten.
    sample_size = min(m, 4 * n)

    if tolerance > 0.0 and requested > 0 and sample_size < m:
        indices = np.random.choice(m, size=sample_size, replace=False)
        indices.sort()

        As = A[indices, :]
        bs = b[indices]
        G = dsyrk(1.0, As.T, trans=0, lower=0)
        sampled_rhs = As.T @ bs

        factor, factor_info = dpotrf(
            G, lower=0, clean=0, overwrite_a=True
        )

        if factor_info == 0:
            x, solve_info = dpotrs(
                factor, sampled_rhs, lower=0, overwrite_b=True
            )

            # Release the large sampled copy before repeatedly scanning A.
            del As, bs, G, sampled_rhs, indices

            if solve_info == 0 and np.all(np.isfinite(x)):
                residual = b - A @ x
                gradient = At @ residual

                z, z_info = dpotrs(
                    factor, gradient, lower=0, overwrite_b=False
                )
                rho = float(np.dot(gradient, z))
                residual_sq = float(np.dot(residual, residual))

                # Equivalent to 1 - 1/(1+eps)^2, evaluated without
                # cancellation for small eps.
                target_fraction = (
                    tolerance * (2.0 + tolerance)
                    / ((1.0 + tolerance) * (1.0 + tolerance))
                )
                certificate_limit = 0.90 * target_fraction

                if (
                    z_info == 0
                    and np.isfinite(rho)
                    and rho >= 0.0
                    and rho <= certificate_limit * residual_sq
                ):
                    return x

                if z_info == 0 and np.isfinite(rho) and rho > 0.0:
                    direction = z.copy()

                    # Bound the speculative work: an unusually poor sample
                    # is cheaper and safer to abandon for the direct fallback.
                    max_iterations = min(requested, 24)

                    for _ in range(max_iterations):
                        Adirection = A @ direction
                        Hdirection = At @ Adirection
                        denominator = float(np.dot(direction, Hdirection))

                        if (
                            not np.isfinite(denominator)
                            or denominator <= 0.0
                        ):
                            break

                        alpha = rho / denominator
                        x += alpha * direction
                        residual -= alpha * Adirection
                        gradient -= alpha * Hdirection

                        z_new, z_info = dpotrs(
                            factor, gradient, lower=0,
                            overwrite_b=False
                        )
                        if z_info != 0:
                            break

                        rho_new = float(np.dot(gradient, z_new))
                        if not np.isfinite(rho_new) or rho_new < 0.0:
                            break

                        residual_sq = float(np.dot(residual, residual))

                        if rho_new <= certificate_limit * residual_sq:
                            # Recompute both quantities from A to eliminate
                            # accumulated CG recurrence error before accepting.
                            true_residual = b - A @ x
                            true_gradient = At @ true_residual
                            true_z, true_info = dpotrs(
                                factor, true_gradient, lower=0,
                                overwrite_b=False
                            )

                            if true_info == 0:
                                true_certificate = float(
                                    np.dot(true_gradient, true_z)
                                )
                                true_residual_sq = float(
                                    np.dot(true_residual, true_residual)
                                )

                                if (
                                    np.isfinite(true_certificate)
                                    and true_certificate >= 0.0
                                    and true_certificate
                                    <= certificate_limit * true_residual_sq
                                ):
                                    return x

                                # Restart CG from the accurately recomputed
                                # gradient if the recurrence was optimistic.
                                if (
                                    np.isfinite(true_certificate)
                                    and true_certificate > 0.0
                                ):
                                    residual = true_residual
                                    gradient = true_gradient
                                    z_new = true_z
                                    rho_new = true_certificate
                                    direction = z_new.copy()
                                    rho = rho_new
                                    continue
                            break

                        if rho_new == 0.0:
                            return x

                        beta = rho_new / rho
                        direction *= beta
                        direction += z_new
                        rho = rho_new

    # Dependable fallback for singular samples, nearly consistent systems,
    # exhausted budgets, or unusually slow preconditioned convergence.
    full_gram = dsyrk(1.0, At, trans=0, lower=0)
    full_rhs = At @ b
    _, x, info = dposv(
        full_gram, full_rhs, lower=0,
        overwrite_a=True, overwrite_b=True
    )

    if info == 0:
        return x

    # Rank-revealing QR preserves the original unregularized objective.
    x, _, _, _ = sla.lstsq(
        A, b,
        cond=None,
        overwrite_a=False,
        overwrite_b=False,
        check_finite=False,
        lapack_driver="gelsy",
    )
    return x
# EVOLVE-BLOCK-END
# ...
\end{lstlisting}

\begin{lstlisting}[
  language=Python,
  caption={Program discovered for the task of size $(30000,2000)$, conditioning $10^{6}$, and accuracy target $10^{-4}$},
  label={code::disco_30k2k_1e6},
  basicstyle=\ttfamily\tiny,  
]
# EVOLVE-BLOCK-START
import numpy as np
from scipy import linalg as sla
from scipy.sparse import csr_matrix

_TINY = 1e-20
_RNG = np.random.default_rng()

def _osnap(m, d):
    """Balanced signed one-sparse embedding with deterministic ||S||_2 <= 1."""
    cols = np.arange(m, dtype=np.int32)
    rows = cols % d
    _RNG.shuffle(rows)
    counts = np.bincount(rows, minlength=d)

    signs = _RNG.integers(0, 2, size=m, dtype=np.int8).astype(np.float64)
    signs *= 2.0
    signs -= 1.0
    signs /= np.sqrt(counts[rows])

    S = csr_matrix((signs, (rows, cols)), shape=(d, m))
    return S, 1.0

def solve(A, b, eps, iter_budget):
    """Solve min_x ||Ax-b|| with OSNAP-preconditioned, certified CGLS."""
    m, n = A.shape
    if n == 0:
        return np.empty(0, dtype=np.float64)

    # A 3n-row sketch provides a highly reliable preconditioner.
    d = min(m, 3 * n)
    S, snorm_bound = _osnap(m, d)
    SA = S @ A
    Sb = S @ b
    del S

    # Compute SA.T @ SA using highly-optimized BLAS GEMM.
    G = SA.T @ SA
    try:
        R = sla.cholesky(G, lower=False, overwrite_a=True, check_finite=False)
    except sla.LinAlgError:
        # Robust fallback if normal-equation rounding loses definiteness.
        R = np.linalg.qr(SA, mode="r")

    if not R.flags.f_contiguous:
        R = np.asfortranarray(R)

    # Cache LAPACK dispatch for the triangular solves outside the loop.
    trtrs = sla.get_lapack_funcs("trtrs", (R,))

    # Warm start: solve R.T @ R @ x = SA.T @ Sb
    rhs = SA.T.dot(Sb)
    z = trtrs(R, rhs, lower=0, trans=1, unitdiag=0, overwrite_b=1)[0]
    x = trtrs(R, z, lower=0, trans=0, unitdiag=0, overwrite_b=1)[0]

    del SA, Sb, G, rhs, z

    r = b - A.dot(x)
    rnorm_sq = float(np.dot(r, r))
    if rnorm_sq <= _TINY:
        return np.asarray(x, dtype=np.float64)

    # Certified stopping criterion based on the sketch norm.
    e = max(float(eps), 0.0)
    eta_sq = e * (2.0 + e) / ((1.0 + e) ** 2) if e > 0.0 else 0.0
    stop_ratio = eta_sq / (1.21 * snorm_bound * snorm_bound)

    # CGLS on B = A R^{-1}, updating x directly.
    g = A.T.dot(r)
    g_res = trtrs(R, g, lower=0, trans=1, unitdiag=0, overwrite_b=1)[0]
    gamma = float(np.dot(g_res, g_res))

    if gamma <= stop_ratio * rnorm_sq:
        return np.asarray(x, dtype=np.float64)

    p = g_res.copy()
    q = np.empty(m, dtype=np.float64)
    A_T_r = np.empty(n, dtype=np.float64)
    direction_x = np.empty(n, dtype=np.float64)

    budget = max(0, int(iter_budget))

    for _ in range(budget):
        np.copyto(direction_x, p)
        dir_res = trtrs(R, direction_x, lower=0, trans=0, unitdiag=0, overwrite_b=1)[0]

        A.dot(dir_res, out=q)
        qq = float(np.dot(q, q))

        if not np.isfinite(qq) or qq <= _TINY:
            break

        alpha = gamma / qq

        # Perform memory-saving, in-place scaling and additions
        dir_res *= alpha
        x += dir_res
        q *= alpha
        r -= q

        rnorm_sq = float(np.dot(r, r))
        if rnorm_sq <= _TINY:
            return np.asarray(x, dtype=np.float64)

        A.T.dot(r, out=A_T_r)
        at_res = trtrs(R, A_T_r, lower=0, trans=1, unitdiag=0, overwrite_b=1)[0]
        gamma_new = float(np.dot(at_res, at_res))

        if gamma_new <= stop_ratio * rnorm_sq:
            return np.asarray(x, dtype=np.float64)

        if not np.isfinite(gamma_new) or gamma <= _TINY:
            break

        p *= gamma_new / gamma
        p += at_res
        gamma = gamma_new

    return np.asarray(x, dtype=np.float64)
# EVOLVE-BLOCK-END
# ...
\end{lstlisting}

\end{document}

%% file: math_commands.tex
\usepackage{amsmath,amsfonts,bm}

\def\eqref#1{equation~\ref{#1}}

\def\1{\bm{1}}

\DeclareMathAlphabet{\mathsfit}{\encodingdefault}{\sfdefault}{m}{sl}
\SetMathAlphabet{\mathsfit}{bold}{\encodingdefault}{\sfdefault}{bx}{n}



%% file: references.bib
@article{romera2024mathematical,
  title={Mathematical discoveries from program search with large language models},
  author={Romera-Paredes, Bernardino and Barekatain, Mohammadamin and Novikov, Alexander and Balog, Matej and Kumar, M. Pawan and Dupont, Emilien and Ruiz, Francisco J.R. and Ellenberg, Jordan S. and Wang, Pengming and Fawzi, Omar and others},
  journal={Nature},
  volume={625},
  number={7995},
  pages={468--475},
  year={2024},
  publisher={Nature Publishing Group UK London}
}

@article{avron2010blendenpik,
  author  = {Avron, Haim and Maymounkov, Petar and Toledo, Sivan},
  title   = {Blendenpik: Supercharging {LAPACK}'s least-squares solver},
  journal = {SIAM Journal on Scientific Computing},
  volume  = {32},
  number  = {3},
  pages   = {1217--1236},
  year    = {2010}
}

@article{murray2023randnla,
  author        = {Murray, Riley and Demmel, James and Mahoney, Michael W. and Erichson, N. Benjamin and Melnichenko, Maksim and Malik, Osman Asif and Grigori, Laura and Luszczek, Piotr and Derezi{\'n}ski, Micha{\l} and Lopes, Miles E. and Liang, Tianyu and Luo, Hengrui and Dongarra, Jack},
  title         = {Randomized Numerical Linear Algebra: A Perspective on the Field with an Eye to Software},
  journal       = {arXiv preprint arXiv:2302.11474},
  year          = {2023},
}

@article{meng2014lsrn,
  author  = {Meng, Xiangrui and Saunders, Michael A. and Mahoney, Michael W.},
  title   = {{LSRN}: A parallel iterative solver for strongly over- or underdetermined systems},
  journal = {SIAM Journal on Scientific Computing},
  volume  = {36},
  number  = {2},
  pages   = {C95--C118},
  year    = {2014}
}

@article{paige1982lsqr,
  title={{LSQR}: An algorithm for sparse linear equations and sparse least squares},
  author={Paige, Christopher C. and Saunders, Michael A.},
  journal={ACM Transactions on Mathematical Software (TOMS)},
  volume={8},
  number={1},
  pages={43--71},
  year={1982},
  publisher={ACM New York, NY, USA}
}

@article{cho2025surrogate,
  title={Surrogate-based autotuning for randomized sketching algorithms in regression problems},
  author={Cho, Younghyun and Demmel, James and Derezi{\'n}ski, Micha{\l} and Li, Haoyun and Luo, Hengrui and Mahoney, Michael and Murray, Riley},
  journal={SIAM Journal on Matrix Analysis and Applications},
  volume={46},
  number={2},
  pages={1247--1279},
  year={2025},
  publisher={SIAM}
}

@article{huh2026adaptive,
  title={Adaptive {LSQR} Preconditioning from One Small Sketch},
  author={Huh, Jung Eun and Cartis, Coralia and Nakatsukasa, Yuji},
  journal={arXiv preprint arXiv:2604.05065},
  year={2026}
}

@article{xiong2026rl4rla,
  title={{RL4RLA}: Teaching {ML} to Discover Randomized Linear Algebra Algorithms Through Curriculum Design and Graph-Based Search},
  author={Xiong, Jinglong and Liu, Xiaotian and Wang, Ruoxin and Liu, Zihang and Zhou, Yefan and Yan, Yujun and Yang, Yaoqing},
  journal={arXiv preprint arXiv:2605.18004},
  year={2026}
}

@article{epperly2026fast,
  title={Fast randomized least-squares solvers can be just as accurate and stable as classical direct solvers},
  author={Epperly, Ethan N. and Meier, Maike and Nakatsukasa, Yuji},
  journal={Communications on Pure and Applied Mathematics},
  volume={79},
  number={2},
  pages={293--339},
  year={2026},
  publisher={Wiley Online Library}
}

@article{cartis2021hashing,
  title={Hashing embeddings of optimal dimension, with applications to linear least squares},
  author={Cartis, Coralia and Fiala, Jan and Shao, Zhen},
  journal={arXiv preprint arXiv:2105.11815},
  year={2021}
}

@article{novikov2025alphaevolve,
  title={{AlphaEvolve}: A coding agent for scientific and algorithmic discovery},
  author={Novikov, Alexander and V{\~u}, Ng{\^a}n and Eisenberger, Marvin and Dupont, Emilien and Huang, Po-Sen and Wagner, Adam Zsolt and Shirobokov, Sergey and Kozlovskii, Borislav and Ruiz, Francisco J.R. and Mehrabian, Abbas and others},
  journal={arXiv preprint arXiv:2506.13131},
  year={2025}
}

@inproceedings{lange2026shinkaevolve,
  title={{ShinkaEvolve}: Towards open-ended and sample-efficient program evolution},
  author={Lange, Robert and Imajuku, Yuki and Cetin, Edoardo},
  booktitle={International Conference on Learning Representations},
  volume={2026},
  pages={74026--74078},
  year={2026}
}

@article{sharma2025openevolve,
  title={{OpenEvolve}: an open-source evolutionary coding agent},
  author={Sharma, Asankhaya},
  journal={GitHub repository https://github.com/algorithmicsuperintelligence/openevolve},
  year={2025}
}

@inproceedings{liu2026skydiscover,
  title={{SkyDiscover}: A flexible, adaptive framework for {AI}-driven scientific and algorithmic discovery},
  author={Liu, Shu and Cemri, Mert and Agarwal, Shubham and Krentsel, Alexander and Naren, Ashwin and Mang, Qiuyang and Li, Zhifei and Gupta, Akshat and Maheswaran, Monishwaran and Cheng, Audrey and others},
  booktitle={Proceedings of the ACM Conference on AI and Agentic Systems},
  pages={1223--1227},
  year={2026}
}

@article{georgiev2025mathematical,
  title={Mathematical exploration and discovery at scale},
  author={Georgiev, Bogdan and G{\'o}mez-Serrano, Javier and Tao, Terence and Wagner, Adam Zsolt},
  journal={arXiv preprint arXiv:2511.02864},
  year={2025}
}

@article{press2026algotune,
  title={{AlgoTune}: Can Language Models Speed Up General-Purpose Numerical Programs?},
  author={Press, Ori and Amos, Brandon and Zhao, Haoyu and Wu, Yikai and Ainsworth, Samuel and Krupke, Dominik and Kidger, Patrick and Sajed, Touqir and Stellato, Bartolomeo and Park, Jisun and others},
  journal={Advances in Neural Information Processing Systems},
  volume={38},
  year={2025}
}

@article{aoudia2026ai,
  title={The {AI} Telco Engineer: Toward Autonomous Discovery of Wireless Communications Algorithms},
  author={A{\"\i}t Aoudia, Fay{\c{c}}al and Hoydis, Jakob and Cammerer, Sebastian and Maggi, Lorenzo and Marti, Gian and Keller, Alexander},
  journal={arXiv preprint arXiv:2604.19803},
  year={2026}
}

@article{golub1965numerical,
  title={Numerical methods for solving linear least squares problems},
  author={Golub, Gene},
  journal={Numerische mathematik},
  volume={7},
  number={3},
  pages={206--216},
  year={1965},
  publisher={Springer}
}

@book{bjorck2024numerical,
  title={Numerical methods for least squares problems},
  author={Bj{\"o}rck, {\AA}ke},
  edition={2},
  year={2024},
  publisher={SIAM}
}

@article{hestenes1952methods,
  title={Methods of conjugate gradients for solving linear systems},
  author={Hestenes, Magnus R. and Stiefel, Eduard},
  journal={Journal of research of the National Bureau of Standards},
  volume={49},
  number={6},
  pages={409--436},
  year={1952}
}

@article{rokhlin2008fast,
  title={A fast randomized algorithm for overdetermined linear least-squares regression},
  author={Rokhlin, Vladimir and Tygert, Mark},
  journal={Proceedings of the National Academy of Sciences},
  volume={105},
  number={36},
  pages={13212--13217},
  year={2008},
  publisher={National Academy of Sciences}
}

@article{martinsson2020randomized,
  title={Randomized numerical linear algebra: Foundations and algorithms},
  author={Martinsson, Per-Gunnar and Tropp, Joel A.},
  journal={Acta Numerica},
  volume={29},
  pages={403--572},
  year={2020},
  publisher={Cambridge University Press}
}

@inproceedings{nelson2013osnap,
  title={{OSNAP}: Faster numerical linear algebra algorithms via sparser subspace embeddings},
  author={Nelson, Jelani and Nguy{\^e}n, Huy L.},
  booktitle={2013 {IEEE} 54th annual symposium on foundations of computer science},
  pages={117--126},
  year={2013},
  organization={IEEE}
}

@misc{susy_279,
  author       = {Whiteson, Daniel},
  title        = {{SUSY}},
  year         = {2014},
  howpublished = {UCI Machine Learning Repository},
  note         = {{DOI}: https://doi.org/10.24432/C54606}
}

@article{baldi2014searching,
  title={Searching for exotic particles in high-energy physics with deep learning},
  author={Baldi, Pierre and Sadowski, Peter and Whiteson, Daniel},
  journal={Nature communications},
  volume={5},
  number={1},
  pages={4308},
  year={2014},
  publisher={Nature Publishing Group UK London}
}

@article{epperly2024fast,
  title={Fast and forward stable randomized algorithms for linear least-squares problems},
  author={Epperly, Ethan N.},
  journal={SIAM Journal on Matrix Analysis and Applications},
  volume={45},
  number={4},
  pages={1782--1804},
  year={2024},
  publisher={SIAM}
}

@article{bergstra2011algorithms,
  title={Algorithms for hyper-parameter optimization},
  author={Bergstra, James and Bardenet, R{\'e}mi and Bengio, Yoshua and K{\'e}gl, Bal{\'a}zs},
  journal={Advances in neural information processing systems},
  volume={24},
  year={2011}
}

@inproceedings{akiba2019optuna,
  title={Optuna: A next-generation hyperparameter optimization framework},
  author={Akiba, Takuya and Sano, Shotaro and Yanase, Toshihiko and Ohta, Takeru and Koyama, Masanori},
  booktitle={Proceedings of the 25th ACM SIGKDD international conference on knowledge discovery \& data mining},
  pages={2623--2631},
  year={2019}
}

@book{golub2013matrix,
  title={Matrix computations},
  author={Golub, Gene H. and Van Loan, Charles F.},
  edition={4},
  year={2013},
  publisher={JHU press}
}

@book{saad2003iterative,
  title={Iterative methods for sparse linear systems},
  author={Saad, Yousef},
  year={2003},
  publisher={SIAM}
}

@article{arioli2015preconditioning,
  title={Preconditioning linear least-squares problems by identifying a basis matrix},
  author={Arioli, Mario and Duff, Iain S.},
  journal={SIAM Journal on Scientific Computing},
  volume={37},
  number={5},
  pages={S544--S561},
  year={2015},
  publisher={SIAM}
}

@article{rudi2017falkon,
  title={{FALKON}: An optimal large scale kernel method},
  author={Rudi, Alessandro and Carratino, Luigi and Rosasco, Lorenzo},
  journal={Advances in neural information processing systems},
  volume={30},
  year={2017}
}

@article{wooldridge2009introductory,
  title={Introductory econometrics: A modern approach},
  author={Wooldridge, Jeffrey M.},
  year={2009},
  publisher={South-Western cengage learning Mason, OH}
}

@book{hastie2009elements,
  title={The elements of statistical learning: data mining, inference, and prediction},
  author={Hastie, Trevor and Tibshirani, Robert and Friedman, Jerome H.},
  volume={2},
  year={2009},
  publisher={Springer}
}

@book{vittinghoff2004regression,
  title={Regression methods in biostatistics: linear, logistic, survival, and repeated measures models},
  author={Vittinghoff, Eric and Glidden, David V. and Shiboski, Stephen C. and McCulloch, Charles E.},
  year={2005},
  publisher={Springer}
}

@book{bevington2003data,
  title={Data reduction and error analysis for the physical sciences},
  author={Bevington, Philip R. and Robinson, D. Keith},
  edition={3},
  year={2003},
  publisher={McGraw-Hill New York}
}

@book{hansen2021computed,
  title={Computed tomography: algorithms, insight, and just enough theory},
  author={Hansen, Per Christian and J{\o}rgensen, Jakob and Lionheart, William R.B. (eds.)},
  year={2021},
  publisher={SIAM}
}

@book{aster2005parameter,
  title={Parameter estimation and inverse problems},
  author={Aster, Richard C. and Borchers, Brian and Thurber, Clifford H.},
  volume={90},
  year={2005},
  publisher={Elsevier academic press Amsterdam}
}

@book{feigelson2012modern,
  title={Modern statistical methods for astronomy: with R applications},
  author={Feigelson, Eric D. and Babu, G. Jogesh},
  year={2012},
  publisher={Cambridge}
}

@article{clarkson2017low,
  title={Low-rank approximation and regression in input sparsity time},
  author={Clarkson, Kenneth L. and Woodruff, David P.},
  journal={Journal of the ACM (JACM)},
  volume={63},
  number={6},
  pages={1--45},
  year={2017},
  publisher={ACM New York, NY, USA}
}

@article{tropp2011improved,
  title={Improved analysis of the subsampled randomized {Hadamard} transform},
  author={Tropp, Joel A.},
  journal={Advances in Adaptive Data Analysis},
  volume={3},
  number={01n02},
  pages={115--126},
  year={2011},
  publisher={World Scientific}
}

@article{williams2009roofline,
  title={Roofline: an insightful visual performance model for multicore architectures},
  author={Williams, Samuel and Waterman, Andrew and Patterson, David},
  journal={Communications of the ACM},
  volume={52},
  number={4},
  pages={65--76},
  year={2009},
  publisher={ACM New York, NY, USA}
}

@article{gholami2024ai,
  title={{AI} and memory wall},
  author={Gholami, Amir and Yao, Zhewei and Kim, Sehoon and Hooper, Coleman and Mahoney, Michael W. and Keutzer, Kurt},
  journal={{IEEE} Micro},
  volume={44},
  number={3},
  pages={33--39},
  year={2024},
  publisher={IEEE}
}

@article{dongarra1990set,
  title={A set of level 3 basic linear algebra subprograms},
  author={Dongarra, Jack J. and Du Croz, Jeremy and Hammarling, Sven and Duff, Iain S.},
  journal={ACM Transactions on Mathematical Software (TOMS)},
  volume={16},
  number={1},
  pages={1--17},
  year={1990},
  publisher={ACM New York, NY, USA}
}

@book{anderson1999lapack,
  title={LAPACK users' guide},
  author={Anderson, Edward and Bai, Zhaojun and Bischof, Christian and Blackford, L. Susan and Demmel, James and Dongarra, Jack and Du Croz, Jeremy and Greenbaum, Anne and Hammarling, Sven and McKenney, Alan and others},
  edition={3},
  year={1999},
  publisher={SIAM}
}

@article{karlson1997estimation,
    author = {Karlson, Rune and W{\'a}lden, Bertil},
    title = {Estimation of optimal backward perturbation bounds for the linear least squares problem},
    journal = {BIT Numerical Mathematics},
    year = {1997},
    pages = {862--869}
}

@InProceedings{ma14statistical,
  title = 	 {A Statistical Perspective on Algorithmic Leveraging},
  author = 	 {Ma, Ping and Mahoney, Michael and Yu, Bin},
  booktitle = 	 {Proceedings of the 31st International Conference on Machine Learning},
  pages = 	 {91--99},
  year = 	 {2014}
}

@inproceedings{eggensperger2013towards,
  title={Towards an empirical foundation for assessing bayesian optimization of hyperparameters},
  author={Eggensperger, Katharina and Feurer, Matthias and Hutter, Frank and Bergstra, James and Snoek, Jasper and Hoos, Holger and Leyton-Brown, Kevin and others},
  booktitle={NIPS workshop on Bayesian Optimization in Theory and Practice},
  volume={10},
  number={3},
  pages={1--5},
  year={2013}
}

@article{greenberg2023optimization,
  title={Optimization or architecture: How to hack kalman filtering},
  author={Greenberg, Ido and Yannay, Netanel and Mannor, Shie},
  journal={Advances in Neural Information Processing Systems},
  volume={36},
  pages={50482--50505},
  year={2023}
}

@inproceedings{saketos2026data,
  title={Data-Driven Discovery of Interpretable Kalman Filter Variants Through Large Language Models and Genetic Programming},
  author={Saketos, Vasileios and Kaltenbach, Sebastian and Litvinov, Sergey and Koumoutsakos, Petros},
  booktitle={International Conference on the Applications of Evolutionary Computation (Part of EvoStar)},
  pages={200--216},
  year={2026},
  organization={Springer}
}
